\documentclass[smallcondensed]{svjour3}
\usepackage{mathrsfs,amsmath,longtable}
\usepackage{graphicx}
\usepackage{mathptmx}
\graphicspath{{../Figs/figures170311/}
              {../Figs/}}
\usepackage{amssymb}
\usepackage{color}

\numberwithin{equation}{section}
\newcommand{\blockmatrix}[3]{%These end of the line comments are neccessary
\begin{minipage}[t][#2][c]{#1}%
\center%
#3%
\end{minipage}%
}%
\newcommand{\fblockmatrix}[3]{%
\fbox{%
\begin{minipage}[t][#2][c]{#1}%
\center%
#3%
\end{minipage}%
}%
}

\newcommand{\ep}[1]{\ensuremath{\varepsilon^{#1}}}

\journalname{Journal ...}

\begin{document}

%\title{ A numerical study of radial basis function methods for solving the Rosenau equation
%\title{Fictitious point and resampling radial basis function methods for solving the Rosenau equation
\title{Radial basis function methods for the Rosenau equation and other higher order PDEs\thanks{Alfa Heryudono was supported by the European Commission CORDIS Marie Curie FP7 program Grant \#235730 and National Science Foundation
DMS Grant \#1552238.}
}

\titlerunning{RBF methods for the Rosenau equation}

\author{Ali Safdari-Vaighani \and Elisabeth Larsson \and Alfa Heryudono }
\authorrunning{A. Safdari-Vaighani, E. Larsson, A. Heryudono}

\institute{
A. Safdari-Vaighani\at
Department of Mathematics,
Allameh Tabataba'i University,
Tehran, Iran\\
\email{asafdari@atu.ac.ir}
\and
E. Larsson
\at
Department of Information Technology,
Uppsala University,
Uppsala, Sweden\\
\email{elisabeth.larsson@it.uu.se}
\and
A. Heryudono\at
Department of Mathematics,
University of Massachusetts Dartmouth,
Dartmouth, Massachusetts, USA\\
\email{aheryudono@umassd.edu}
}

\date{Received: date / Accepted: date}

\maketitle

%%%%%%%%%%%%%%%%%%%%%%%%%%%%%%%%%%%%%%%%%%%%%%%%%%%%%%%%%%%%%%%%%%%%%%%%%%%%%%%%%%%%%%%%%%%%%%%%%%%%%%%%%
\begin{abstract}
Meshfree methods based on radial basis function (RBF) approximation are of interest for numerical solution of partial differential equations (PDEs) because they are flexible with respect to the geometry of the computational domain, they can provide high order convergence, they are not more complicated for problems with many space dimensions and they allow for local refinement. The aim of this paper is to show that the solution of the Rosenau equation, as an example of an initial-boundary value problem with multiple boundary conditions, can be implemented using RBF approximation methods. We extend the fictitious point method and the resampling method to work in combination with an RBF collocation method. Both approaches are implemented in one and two space dimensions. The accuracy of the RBF fictitious point method is analysed partly theoretically and partly numerically. The error estimates indicate that a high order of convergence can be achieved for the Rosenau equation. The numerical experiments show that both methods perform well. In the one-dimensional case, the accuracy of the RBF approaches is compared with that of a pseudospectral resampling method, showing similar or slightly better accuracy for the RBF methods. In the two-dimensional case, the Rosenau problem is solved both in a square domain and in a starfish-shaped domain, to illustrate the capability of the RBF-based methods to handle irregular geometries.

\end{abstract}
% 35G31 Initial-boundary value problems for nonlinear higher-order equations

%
\section{Introduction}
The Rosenau equation has become an established research subject in the field of mathematical
physics since its introduction in the late 80s by Philip Rosenau~\cite{rosenau}.
The equation is intended to overcome shortcomings of the already
famous Korteweg--de Vries (KdV) equation~\cite{forwhi} in describing phenomena of
solitary wave interaction. Knowledge about this interaction,
particularly when two or more wave packets called \emph{solitons}
are colliding with one another, is indispensable in digital
transmission through optical fibers. As data carriers, we need
solitons that interact ``cleanly'' in the sense that none of the solitons
loose any information, shape, or other conserved
quantities, when they pass through each other. One may consult \cite{cipra}
for a fascinating history behind this subject.
The Rosenau equation in its general form is given by
\begin{align}
u_t+\alpha(\mathbf{x},t)\Delta^2u_{t}=\nabla\cdot
g(u),\quad(\mathbf{x},t)\in \Omega\times(0,\,T], \label{eq:Rosenau}
\end{align}
where $\Omega$ is a bounded domain in $\mathbb{R}^d$ ($d\leq 3$), the coefficient $\alpha(x,t)\geq\alpha_0>0$ is bounded in the domain $\Omega\times[0,T]$, and the nonlinear function $g(u)$ is polynomial of degree $q+1$, $q\geq 1$.
Multiple boundary conditions are required at the boundary $\partial \Omega$, such as
\begin{align}
u(\mathbf{x},t)&=f_1(\mathbf{x},t),\quad (\mathbf{x},t)\in
\partial\Omega\times(0,\,T],\label{eq:bc1}\\
\frac{\partial u}{\partial n}(\mathbf{x},t)&=f_2(\mathbf{x},t),\quad(\mathbf{x},t)\in
\partial\Omega\times(0,\,T],\label{eq:bc2}
\end{align}
where ${n}$ is the outward normal direction from the boundary, and we need an initial condition
\begin{align*}
u(\mathbf{x},0)=f_0(\mathbf{x}),\quad x\in \Omega. \label{eq:ic}
\end{align*}

The well-posedness of the Rosenau equation has been studied theoretically by Park~\cite{park,phdpark}. Yet in practice, the equation poses difficulties to solve numerically due to the presence of non-linearity, high spatial derivatives, multiple boundary conditions, and mixed time and space derivatives. Numerical studies based on Galerkin formulations can be found in~\cite{ChuHa94,kimlee,leeahn,ChuPa01}, and numerical studies based on finite difference methods are found in~\cite{Chung98,OAAK08,HuZhe10,AtOm15}.

%\comment{the diffusion coefficient $\beta$
%must be negative in order that initial value problem be well-posed in forwards time.}
% The purpose of this paper.
% RBF for the Rosenau equation
% More than 1-D
% General problems with mutliple bcs

The objective of this paper is to derive numerical methods based on radial basis function (RBF) collocation methods~\cite{LaFo03,FoFly15} for the Rosenau equation, that can be applied to problems in one, two, and three space dimensions, for non-trivial geometries. These methods will also be applicable to other higher order partial differential equations. We derive and implement a fictitious point RBF collocation method and a resampling RBF collocation method, and perform experiments in one and two space dimensions. We investigate the accuracy and behavior of the derived methods theoretically and numerically. We also compare the RBF methods with a pseudospectral (PS) method~\cite{Fornberg96,Tref00} with respect to accuracy in one space dimension.

The outline of the paper is as follows: In Section~\ref{sec:rbf}, a basic RBF collocation scheme is introduced. Section~\ref{sec:bc} describes different approaches to handle the multiple boundary conditions. Then in Section~\ref{sec:theory} the approximation properties of the RBF method for the one-dimensional Rosenau problem are analyzed theoretically. Implementation aspects are discussed in Section~\ref{sec:matlab}, followed by numerical results in Section~\ref{sec:num}. Finally, Section~\ref{sec:conc} contains conclusions and discussion.

\section{The basic RBF collocation scheme}
\label{sec:rbf}
The approaches for handling multiple boundary conditions implemented in this paper are combined with an RBF collocation method. In this section, we introduce the general notation and quantities we need for RBF approximation of (time-dependent) PDEs. We start from given scalar function values $u(x_j)$ at scattered distinct node locations $x_j\in\Bbb{R}^d$, $j=1,\ldots,N$. We assume that the function is approximated by a standard RBF interpolant
\begin{equation}
s(x)=\sum^{N}_{j=1}\lambda_j\phi(\|x-x_j\|),
\label{eq:rbfint}
\end{equation}
where $\|\cdot\|$ is the Euclidean norm,  $\phi$ is a real-valued function such as
the inverse multiquadric $\phi(r)=\frac{1}{\sqrt{\varepsilon^2r^2+1}}$.
The coefficients $\lambda_j\in \Bbb{R}$
are determined by
the interpolation conditions $s(x_j)=u(x_j),$ $j=1,\ldots,N$. On matrix form we have
\begin{equation}
A\bar{\lambda}=\bar{u}, \label{eq:sys}
\end{equation}
where the matrix elements $A_{ij}=\phi(\|x_i-x_j\|),~i,j=1,\ldots,N$, the vector $\bar{\lambda}=[\lambda_1,\ldots ,\lambda_N]^T$, and
$\bar{u}=[u(x_1),\ldots, u(x_N)]^T$. When solving a PDE, we prefer to work with the discrete function values instead of the coefficients. Using~\eqref{eq:rbfint} and~\eqref{eq:sys} together, we see that the interpolant can be written
\begin{equation}
s(x)=\bar{\phi}(x)\bar{\lambda}=\bar{\phi}(x)A^{-1}\bar{u},
\label{eq:phi}
\end{equation}
where $\bar{\phi}(x)=[\phi(\|x-x_1\|), \ldots, \phi(\|x-x_N\|)]$, assuming that $A$ is non-singular. This holds for commonly used
RBFs such as Gaussians, inverse multiquadrics and multiquadrics~\cite{Scho38,micch} for distinct node points $x_j$. We can furthermore, use~\eqref{eq:phi} to identify cardinal basis functions such that we can write the approximation on the FEM like form
\begin{equation}
s(x)=\bar{\psi}(x)\bar{u}=\sum^N_{j=1}\psi_j(x)u(x_j),
\label{eq:phipsi}
\end{equation}
where $\bar{\phi}(x)A^{-1}=\bar{\psi}(x)=[\psi_1(x),\ldots,\psi_N(x)]$. Because our final target is to solve PDEs, we need to look at how to apply a linear operator $\mathcal{L}$ to the RBF approximation to compute $s_{\mathcal{L}}=[\mathcal{L}s(x_1),\ldots, \mathcal{L}s(x_N)]^T$. In cardinal form, we get
\begin{equation}
\mathcal{L}{s}(x)=\mathcal{L}\bar{\psi}(x)\bar{u} =\sum^N_{j=1}\mathcal{L}\psi_j(x)u(x_j).
\end{equation}
Using relation~(\ref{eq:phi}), the differentiation matrix $\Psi_\mathcal{L}=\{\mathcal{L}\psi_j(x_i)\}_{i,j=1,\ldots,N}$ under operator $\mathcal{L}$ is given by 
\begin{equation}
\Psi_\mathcal{L}=\Phi_\mathcal{L}A^{-1},
\label{eq:diffmat}
\end{equation}
where $\Phi_\mathcal{L}=\{\mathcal{L}\phi(\|x_i-x_j\|)\}_{i,j=1,\ldots,N}$.

For time-dependent PDE problems, we use the RBF approximation in space and then discretize the time interval. The solution $u(x,t)$ is approximated by
\begin{equation}
s(x,t)=\sum^N_{j=1}\psi_j(x)u_j(t), \label{eq:rbfinterp}
\end{equation}
where $u_j(t)\approx u(x_j,t)$ are the unknown functions to determine.

\section{Dealing with multiple boundary conditions}
\label{sec:bc}
If we consider the one-dimensional version of equations~(\ref{eq:Rosenau})--(\ref{eq:bc2}) on a symmetric interval $x\in[-L,\,L]$ we have
\begin{equation}
u_t+\alpha(x,t)u_{xxxxt}=
g_u(u)u_x,\quad(x,t)\in [-L,\,L]\times(0,\,T],
\label{eq:Rosenau1D}
\end{equation}
with boundary conditions
\begin{align}
u(\pm L,t)&=f_1(\pm L,t),\quad t\in (0,\,T],\\
u_x(\pm L,t)&=f_2(\pm L,t),\quad t\in (0,\,T].\label{eq:1Dic}
\end{align}
%and initial condition
%\begin{equation}
%u(x,0)=f_0(x),\quad x\in \Omega. \label{eq:ic}
%\end{equation}
%
Even for the one-dimensional case, how to implement multiple boundary
conditions for a time-dependent global collocation problem is not
obvious. In our case, we need to enforce two boundary conditions
at each end point resulting in a total of four boundary conditions at the two
boundary points. Collocating the PDE at all interior node points leads to a situation where we have more
equations than node points. That is, unless we accept an overdetermined system, we need to either increase the number of degrees of freedom or discard some of the equations.

In fact, the subtleties of implementing multiple BCs are not tied
to RBF methods. They have been actively researched in conjunction with other
global collocations methods, particularly pseudospectral methods,
since the 70s. We list the following five popular methods:
\begin{enumerate}
\item Mixed hard and weak BCs~\cite{hesgot}
\item Spectral penalty method~\cite{hesta}
\item Transforming to lower order system~\cite{singer}
\item Fictitious/ghost point method~\cite{fornberg}
\item Resampling method~\cite{Drishal}
\end{enumerate}
In this paper, we only consider methods (3)--(5) as we currently do not have a way to find penalty parameters for methods (1) and (2)
that give numerically stable solutions.

\subsection{Transforming to lower order system}
A common approach when solving PDEs
containing high order derivatives is to transform them into a system with
lower order derivatives. By letting $w = u_x$, the Rosenau equation~(\ref{eq:Rosenau1D}) is
transformed into
\begin{align*}
u_t + \alpha(x,t)w_{xxxt}&= wg_u(u)\\
w_t - u_{xt} &= 0
\end{align*}
with boundary conditions $u(\pm L,t)=f_1(\pm L,t)$
and $w(\pm L,t)=f_2(\pm L,t)$. The advantage of this method
is that the Neumann conditions for $u$ at the boundaries become
Dirichlet conditions for $w$. However, the system to solve becomes
twice as large, as we need a total of $2N$ degrees of freedom for $u$
and $w$. Due to this reason, and especially for global RBFs where
differentiation matrices are dense, we are not pursuing this
method. However, for RBF methods in finite difference mode where
differentiation matrices are sparse, this method may still be worth
trying.

\subsection{Fictitious point method}
\label{sec:fp}
Fictitious or ghost point methods have been commonly used as a way to
enforce multiple boundary conditions in finite difference methods.
The implementation for global collocation methods such as
pseudospectral methods is due to Fornberg~\cite{fornberg}.

The Dirichlet conditions~(\ref{eq:bc1}) can be imposed by
fixing the values for
$u(x_1)$ and $u(x_N)$, but for the Neumann conditions~(\ref{eq:bc2}) we use the fictitious point approach
proposed by Fornberg~\cite {fornberg}, and introduce two additional
points at some arbitrary locations such as $x_0$ and $x_{N+1}$.

We introduce an RBF approximation $s(x,t)$ as
in~(\ref{eq:rbfinterp}), extended to include the fictitious
points, for the spatial approximation of the solution $u(x,t)$,

\begin{equation}
s(x,t)=\sum^{N+1}_{j=0}\psi_j(x)u_j(t). \label{eq:spacedisc}
\end{equation}

Loosely following the fictitious point approach, we will modify this ansatz so that the boundary conditions are
fulfilled. Conditions~(\ref{eq:bc1}) are easily fulfilled by replacing $u_j(t)$ with $f_1(x_j,t)$ for $j=1,N$.
For the conditions~(\ref{eq:bc2}), we need to formally solve a linear system. Define the vectors $S_f=[u_0(t),\, u_{N+1}(t)]^T$
with values at the two fictitious points, and $S_d=[u_2(t),\ldots, u_{N-1}(t)]^T$ containing the approximate solution
values at points in the interior of the domain, then we have

\begin{equation}
\underbrace{
\left(
\begin{array}{cc}
\psi_0^\prime(x_1) &
\psi_{N+1}^\prime(x_1)\\
\psi_0^\prime(x_N) &
\psi_{N+1}^\prime(x_N)
\end{array}
\right)
}_{B_f}
S_f+
\underbrace{
\left(
\begin{array}{ccc}
\psi_2^\prime(x_1) & \cdots &
\psi_{N-1}^\prime(x_1)\\
\psi_2^\prime(x_N) & \cdots &
\psi_{N-1}^\prime(x_N)
\end{array}
\right)
}_{B_d}
S_d
+\underbrace{
\left(
\begin{array}{cc}
\psi_1^\prime(x_1) &
\psi_N^\prime(x_1)\\
\psi_1^\prime(x_N) &
\psi_N^\prime(x_N)\\
\end{array}
\right)}_{B_b}
F_1(t)
=
F_2(t),
\label{eq:Uf}
\end{equation}
where $F_j(t)=[f_j(x_1,t),\,f_j(x_N,t)]^T$. Inserting the boundary
values $F_1(t)$ and the expression we get for $S_f$ by solving~\eqref{eq:Uf} into~(\ref{eq:spacedisc}) leads to
\begin{align}
s(x,t)&= \left([\psi_2(x),\ldots,\psi_{N-1}(x)]-[\psi_0(x),\, \psi_{N+1}(x)]\,B_f^{-1}B_d\right)\,S_d \nonumber\\
&+\left([\psi_1(x),\, \psi_{N}(x)]-[\psi_0(x),\,
\psi_{N+1}(x)]\,B_f^{-1}B_b \right)F_1(t)+ [\psi_0(x),\,
\psi_{N+1}(x)]\,B_f^{-1}F_2(t). \label{eq:modint1}
\end{align}
This expression is awkward to work with directly. We introduce the shorthand notation

\begin{equation}
s(x,t)= \sum^{N-1}_{j=2}\tilde{\psi}_j(x)u_j(t)+F(x,t),
\label{eq:modint2}
\end{equation}
where $\tilde{\psi}_j(x)$ and $F(x,t)$ can be directly identified from~(\ref{eq:modint1}).
In this simple two point boundary case, we can actually derive the explicit form of the modified basis for illustration. This yields

\begin{align}
\tilde{\psi}_j(x)=\psi_j(x)&- \frac{\psi_{N+1}^\prime(x_N)
\psi_{j}^\prime(x_1)-
\psi_{N+1}^\prime(x_1)
\psi_{j}^\prime(x_N)
}
{
\psi_{N+1}^\prime(x_N)\psi_0^\prime(x_1)-
\psi_{N+1}^\prime(x_1)
\psi_0^\prime(x_N)
}\psi_0(x)\nonumber\\[3pt]
&+ \frac{ \psi_{0}^\prime(x_N) \psi_{j}^\prime(x_1)-
\psi_{0}^\prime(x_1) \psi_{j}^\prime(x_N) } {
\psi_{N+1}^\prime(x_N) \psi_0^\prime(x_1) - \psi_{N+1}^\prime(x_1)
\psi_0^\prime(x_N) } \psi_{N+1}(x).
\end{align}

In order to use the RBF approximation~(\ref{eq:modint2}) for a PDE problem, we need to compute the effect of applying a spatial differential operator $\mathcal{L}$ at the interior node points. That is, we need a method to evaluate $\mathcal{L}\tilde{\psi}_j(x_i)$, $i,j=2,\ldots,N-1$, and $\mathcal{L}F(x_i,t)$, $i=2,\ldots,N-1$.
This is done in two steps. First, we use~(\ref{eq:diffmat}) to compute $\Psi_\mathcal{L}$ for interior node points $x_i$, $i=2,\ldots,N-1$. Note however that we include all basis functions $\psi_j(x)$, $j=0,\ldots,N+1$. Then we extract the columns pertaining to the fictitious points into $\Psi_{\mathcal{L},f}$, the columns pertaining to the boundary points into $\Psi_{\mathcal{L},b}$, and the remaining columns into $\Psi_{\mathcal{L},d}$. Then the modified differentiation matrix and the contribution in the forcing function can be computed as

\begin{equation}
\tilde{\Psi}_{\mathcal{L}}=\Psi_{\mathcal{L},d}-\Psi_{\mathcal{L},f}B_f^{-1}B_d,
\label{eq:moddiff}
\end{equation}

\begin{equation}
[F_{\mathcal{L}}(x_2,t),\ldots,F_{\mathcal{L}}(x_{N-1},t)]^T=(\Psi_{\mathcal{L},b}-\Psi_{\mathcal{L},f}B_f^{-1}B_b)F_1(t)+\Psi_{\mathcal{L},f}B_f^{-1}F_2(t).
\label{eq:bondim}
\end{equation}
Note that from~(\ref{eq:diffmat}), if no operator is applied, we have $\Psi_d=I$ and $\Psi_f=\Psi_b=0$ leading to $F(x_i,t)=F_t(x_i,t)=0$.

Collocating the modified RBF approximation~(\ref{eq:modint2}) with
the PDE~(\ref{eq:Rosenau}) at the node points leads to the system
of ODEs
\begin{align}
u_i^\prime(t)&+\sum_{j=2}^{N-1}\alpha(x_i,t)\frac{d^4
\tilde{\psi}_j}{d
x^4}(x_i)u_j^\prime(t)=\sum_{j=2}^{N-1}g_u(u_i(t))\frac{d \tilde{\psi}_j}{d
x}(x_i)u_j(t) \nonumber \\[3pt]
&-\alpha(x_i,t) F_{xxxxt}(x_i,t)+g_u(u_i(t))F_x(x_i,t), \quad i=2,\ldots,N-1.
\end{align}
In matrix form, we get the method of lines formulation

\begin{equation}
\underbrace{(I+A_\alpha(t)\tilde{\Psi}_{xxxx})}_{Q(t)}S_d^\prime=\underbrace{G_u(S_d)\tilde{\Psi}_x}_{D(S_d)}S_d+\underbrace{G_u(S_d)F_x(t) - A_\alpha(t)F_{xxxx}^\prime(t)}_{F(S_d,t)},
\label{eq:fpmol}
\end{equation}
where the diagonal coefficient matrices are
\begin{eqnarray*}
A_\alpha(t) &=&\mathrm{diag}(\alpha(x_2,t),\ldots,\alpha(x_{N-1},t)),\\
G_u(S_d) &=&\mathrm{diag}(g_u(u_2(t)),\ldots,g_u(u_{N-1}(t))),
\end{eqnarray*}
and the vectors in the right hand side are defined as $F_\mathcal{L}(t)=[F_\mathcal{L}(x_2,t),\ldots,F_\mathcal{L}(x_{N-1},t)]^T$.
The problem~\eqref{eq:fpmol} can be solved by employing a solution method for nonlinear ODE systems.

The coefficient matrix $Q(t)$ is in general invertible but
non-singularity cannot be guaranteed.
Kansa \cite {kansa} argued that if the centers
of the RBFs are distinct and the PDE problem is well-posed, the
coefficient matrix is generally found to be non-singular.
Hon and
Schaback \cite {honscha} showed that occurrences of singular
coefficient matrix are very rare, but do exist.

When $\alpha(x,t)$ is constant, %(for example in the Rosenau--Burgers equation),
the coefficient matrix is constant over time.
Then we can LU-factorize the
coefficient matrix once and use this factorization throughout the time stepping
algorithm.

An alternative to the derivation above is to use the original cardinal basis functions, and include the boundary condition equations in the final system.
Define rectangular identity matrices $I_k$ such that $I_k(S_d,\,S_b,\,S_f)^T=S_k$, for $k=d,b,f$. Also, we order the columns in differentiation matrices in accordance with the order of the unknowns such that $\Psi_\mathcal{L}=[\Psi_{\mathcal{L},d},\,\Psi_{\mathcal{L},b},\Psi_{\mathcal{L},f}]$. Then we can express~\eqref{eq:fpmol} as
\begin{equation}
\left(\begin{array}{c}
I_d+A_{\alpha}\Psi_{xxxx} \\
0\\
0
\end{array}\right)
\left(\begin{array}{c}
S_d^\prime\\
S_b^\prime\\
S_f^\prime
\end{array}\right)=
\left(\begin{array}{c}
G_u(S_d)\Psi_x\\
I_b\\
B_d\ \ B_b\ \ B_f
\end{array}\right)
\left(\begin{array}{c}
S_d\\
S_b\\
S_f
\end{array}\right)-
\left(\begin{array}{c}
0\\
F_1(t)\\F_2(t)
\end{array}\right).
\label{eq:fict}
\end{equation}
%
%%%%%%%%%%% START HERE WITH HARMONIZING THE TWO PARTS %%%%%%%%%%%%
% Use M+A(t). Another U since all points are there. Order? Compare implementation section. Change to have same notation at both places. Implementation section is very specific. Where do we choose alpha? This is not described as resampling. Why don't we do it like that?
\subsection{Resampling method}\label{sec:res}
In the resampling method, we do not add any points as for the fictitious point method of the previous section. The four boundary conditions
are still enforced at the boundary points as algebraic equations, but
the PDE is instead collocated at $N-4$ auxiliary interior points.
%The PDE part is discretized by the RBF interpolant and differentiation matrices at the auxiliary set of point that is different and transformed by RBF center points.
%
Let the solution $u(x,t)$ be approximated in Lagrange form by
\begin{align}
s(x,t)=\sum^{N}_{j=1}\psi_j(x)u_j(t),
\label{eq:rbfinterp2}
\end{align}
where $x_1$ and $x_N$ are boundary points and $x_2,\ldots,x_{N-1}$ are interior points.
The four algebraic equations arising from the boundary conditions are
\begin{align}
u_1(t) = f_1(x_1,t), & \quad
u_N(t) = f_1(x_N,t),\label{eq:algebraic1}\\
\sum_{j=1}^{N}\frac{d\psi_j}{dx}({x}_1)u_j(t) = f_2(x_1,t), & \quad
\sum_{j=1}^{N}\frac{d\psi_j}{dx}({x}_N)u_j(t) = f_2(x_N,t).
\label{eq:algebraic2}
\end{align}
To write the equations on matrix form, we again divide the unknown functions into parts, $S_d$ for interior node points and $S_b$ for boundary node points. Then the boundary conditions can be expressed as
\begin{equation}
\left(\begin{array}{cc}
0 & I_b\\
\tilde{B}_d & \tilde{B}_b
\end{array}\right)
\left(\begin{array}{c}
S_d\\
S_b\\
\end{array}\right)
=
\left(\begin{array}{c}
F_1(t)\\
F_2(t)
\end{array}\right),
\label{eq:resamplingbc}
\end{equation}
where the boundary condition matrices $\tilde{B}_d$ and $\tilde{B}_b$ are analogous to those in the previous section, but with slightly different basis functions as there are no fictitious points.
We have $N$ unknown functions, and we have used four equations for the boundary conditions. This means that we need $N-4$ additional equations. The idea in the resampling method is to collocate the PDE at $N-4$ auxiliary interior points, instead of collocating at the node points.
We define the auxiliary points $\tilde{x}_i$, $i=1,\ldots,N-4$ and collocate the PDE to get the equations
\begin{align}
\sum_{j=1}^{N}\left[\psi_j(\tilde{x}_i)+\alpha(\tilde{x}_i,t)\frac{d^4 \psi_j}{d x^4}(\tilde{x}_i)\right]u_j^\prime(t)=\sum_{j=1}^{N}\left[g_u\left(\sum_{k=1}^{N}\psi_k(\tilde{x}_i)u_k(t)\right)\frac{d \psi_j}{d x}(\tilde{x}_i)\right]u_j(t), %\quad i=1,\ldots,N-4.
\label{eq:differential}
\end{align}
%
% Curly bracket, index
Define the resampling matrices $\Psi_\mathcal{L}^R=\{\mathcal{L}\psi_j(\tilde{x}_i)\}_{i=1,\ldots,N-4,\, j=2,\ldots,N-1,1,N}$ (columns ordered according to the unknowns).
The resampled equations~\eqref{eq:differential}, together with the algebraic equations~\eqref{eq:resamplingbc}, yield an $N \times N$ differential algebraic system of equations (DAE) of index~1,
% Have not explained the S, which is the whole vector
\begin{equation}
\left(\begin{array}{c}
\Psi^R+\tilde{A}_{\alpha}\Psi^R_{xxxx}\\
0\\
0
\end{array}\right)
\left(\begin{array}{c}
S_d^\prime\\
S_b^\prime\\
\end{array}\right)=
\left(\begin{array}{cc}
\multicolumn{2}{c}{G_u(\tilde{S})\Psi^R_x}\\
0 & I_b\\
\tilde{B}_d & \tilde{B}_b
\end{array}\right)
\left(\begin{array}{c}
S_d\\
S_b\\
\end{array}\right)-
\left(\begin{array}{c}
0\\
F_1(t)\\
F_2(t)
\end{array}\right),
\label{eq:resampling}
\end{equation}
where $\tilde{A}_{\alpha}=\mathrm{diag}(\alpha(\tilde{x}_1,t),\ldots,\alpha(\tilde{x}_{N-4},t) )$ and $\tilde{S}=\Psi^R\left(\begin{array}{c}
S_d\\
S_b\\
\end{array}\right)$.

%The $4\times N$ matrix $C$ of the algebraic part has elements $C_{1j}=\psi_j(x_1)=\delta_{1j}$, $C_{2j}=\psi_j(x_N)=\delta_{Nj}$,
%$C_{3j}=\frac{d\psi_j}{dx}({x}_1)$, and $C_{4j}=\frac{d\psi_j}{dx}({x}_N)$. The vector $U=[u_1(t),\ldots,u_N(t)]^T$ contains the approximate solution values at all the node points, and $F(t)=[f_1(x_1),\,f_1(x_N),\,f_2(x_1),\,f_2(x_N)]^T$.

The system of equations \eqref{eq:resampling} can be solved using a differential algebraic solver. See \textsc{DASPK}~\cite{brohinpet,petzold}. An example of how this can be implemented in \textsc{MATLAB} is given in section~\ref{sec:matlab}.

%%%%%%%%%%%%%%%%%%%%%%%%%%%%%%%%%
\subsection{Generalization to more space dimensions}
% Main difference is that hte number of boundary points are more
The main differences when moving to more than one space dimensions is that we have a boundary curve or a boundary surface that is discretized in the same way as the interior of the domain instead of just two boundary points. The formulation of the two methods is in all essential parts the same, and the formulations~\eqref{eq:fict} and~\eqref{eq:resampling} are valid in the same form, but when we before had two boundary points and two fictitious points, we instead have $N_b$ boundary points and $N_b$ fictitious points. Similarly, for the resampling method, we have $2N_b$ boundary conditions, and therefore we collocate the PDE at $N-2N_b$ auxiliary points. Experiments for problems in two space dimensions are presented in Section~\ref{sec:num}.

\section{Error estimates}\label{sec:theory}
% Is also F_{xxxxt}=0? check the dependencies
% How about having d/dt(E+epsilon)? We know how to bound interpolation error
% Then the interpolation error comes out nicely. Lets try it.
% F(x,t)=A(x)F1(t)+B(x)F2(t)
% F_xxxx(x,t)=A_xxxxF1(t)+B_xxxxF2(t) \neq 0
% Integral from zero because E(0)=0.

% Can we say about the eigenvalues of A_xxxx?
% If A has 1/(1+a). Probably this could be tried
%-------------------
% Start with what we are going to do
% Then make a new section ODE system for the discrete approximation error
In this section, we derive semi-discrete error estimates for the one-dimensional problem. The analysis is carried out for the fictitious point approach. For the analysis, we assume that $\alpha>0$ is constant and that the function $g$ is a polynomial of degree $q+1$
with $q\geq 1$.
\subsection{ODE system for the semi-discrete approximation error}
% I will want to define the vectors, perhaps one at a time.
% U_d or S_d, Now because I want to have derivatives.
% First define U_L
% From Rosenau we know that
We work with spatially discrete solution and approximation values evaluated at the interior node points $x_i$, $i=2,\ldots,N-1$.
We denote the exact solution vector and its derivatives by $U_\mathcal{L}(t)=(\mathcal{L}u(x_2,t),\ldots,\mathcal{L}u(x_{N-1},t))^T$.
From the Rosenau equation~\eqref{eq:Rosenau1D}, we have
\begin{equation}
U^\prime+\alpha U_{xxxx}^\prime=G_u(U)U_x,
\label{eq:uexmol}
\end{equation}
For the RBF approximation, we use~\eqref{eq:fpmol}, but to simplify notation we replace $S_d$ with $S$ and the matrix $A_\alpha$ with the constant $\alpha$.
% Simplify the notation for S, let A_a=a
\begin{equation}
(I+\alpha\tilde{\Psi}_{xxxx})S^\prime+\alpha F_{xxxx}^\prime=G_u(S)\tilde{\Psi}_{x}S+G_u(S)F_x.
\label{eq:fpmol2}
\end{equation}
% Now we will introduce an auxiliary Z and describe its derivatives
% Then we will rewrite (4.1) with addition and subtraction
% Perhaps we also introduce epsilon for those terms already here.
We introduce the auxiliary function $z(x,t)$, which interpolates the exact solution at each time and also satisfies the boundary conditions~\eqref{eq:bc1} and~\eqref{eq:bc2},
\begin{equation}
z(x,t)=\sum_{j=2}^{N-1}\tilde{\psi}_j(x)u(x_j,t)+F(x,t).\label{eq:mol1}
\end{equation}
We denote the auxiliary function vector by $Z(t)$ and note that
\begin{equation}
Z_\mathcal{L}=\tilde{\Psi}_\mathcal{L}U+F_\mathcal{L}
\label{eq:aux}
\end{equation}
Furthermore, $Z(t)=U(t)$, $Z'(t)=U'(t)$, and $Z(0)=S(0)=U(0)$, while $Z_\mathcal{L}\neq U_\mathcal{L}$.

%\begin{eqnarray}
%U^\prime+\alpha (U_{xxxx}^\prime-Z_{xxxx}^\prime+Z_{xxxx}^\prime) &=& G_u(U)(U_x-Z_x+Z_x),
%\end{eqnarray}
We define the discrete interpolation error vector $\epsilon(t)=U(t)-Z(t)$. The interpolation error itself is zero at the node points, but its derivatives $\epsilon_\mathcal{L}$ are non-zero. By noting that $U(t) = \epsilon(t) + Z(t)$ and by using~\eqref{eq:aux} for the derivatives of $Z(t)$, we can rewrite~\eqref{eq:uexmol} as
\begin{eqnarray}
U^\prime+\alpha (\epsilon_{xxxx}^\prime+\tilde{\Psi}_{xxxx}U^\prime+F_{xxxx}^\prime) &=&
G_u(U)(\epsilon_{x}+\tilde{\Psi}_{x}U+F_x).
\label{eq:uexmod}
\end{eqnarray}
Finally, we introduce the discrete error $E(t)=U(t)-S(t)$, and subtract~\eqref{eq:fpmol2} from~\eqref{eq:uexmod} to get
% Need to define the error
\begin{equation}
\underbrace{(I+\alpha\tilde{\Psi}_{xxxx})}_{Q}E^\prime + \alpha\epsilon_{xxxx}^\prime =
G_u(U)\epsilon_{x}+ G_u(U)\tilde{\Psi}_{x}U -G_u(S)\tilde{\Psi}_{x}S + (G_u(U)-G_u(S))F_x.
\end{equation}
%% Missing some F-terms here
This equation can be seen as an ODE-system for the error $E(t)$ by writing it as
\begin{equation}
QE^{\prime}(t)=H(t)\label{error1},
\end{equation}
%\begin{equation}
%(I+\alpha(t)\Psi_{xxxx})(U-S)^{\prime}(t)=-\beta\Psi_{xx}(U-S)(t)+H(t)\label{error1},
%\end{equation}
where
$H(t)=-\alpha\epsilon^{\prime}_{xxxx}(t)+
G_u(U)\tilde{\Psi}_xU-G_u(S)\tilde{\Psi}_xS+G_u(U)\epsilon_x(t)+
(G_u(U)-G_u(S))F_x(t)$. In the following, we will consider $H(t)$ as a forcing function.
%
%$g_u(U(t))\epsilon_x(t)+g_u(U(t))\Psi_xU(t)-g_u(S(t))\Psi_xS(t)
%+F_x(t)[g_u(U(t))-g_u(S(t))]$,
%$F(t)=[F(x_2,t),\ldots,F(x_{N-1},t)]^T$ and
%$\epsilon(t)=U(t)-Z(t)$.
%
For a discussion of the non-singularity of $Q$, see Section~\ref{sec:fp}. The system of ODEs~(\ref{error1}) can be formally integrated to yield
%condition~(\ref{ini:mol}) and boundedness of $\alpha(t)$ is
%\begin{equation}
%(U-S)(t)\leq\int_{0}^{t}\exp^{P(t-\tau)}Q
%H(\tau)d\tau\label{ode:ersol},
%\end{equation}
\begin{equation}
E(t)=\int_{0}^{t}Q^{-1}H(\tau)d\tau.
\label{ode:ersol}
\end{equation}
In the following subsections, we will look into each part of the forcing function.
\subsection{Estimates for the non-linear term}
In order to determine the influence of the non-linear term on convergence and stability, we consider the particular form $g(u)=u^{q+1}$, $q\geq 1$. This matches the functions typically used in the literature. Furthermore, we will use an estimate by Park from \cite{park} showing that with this form of $g(u)$,
\begin{equation}
|u(x,t)|\leq
C(1+t)^{1-\frac{q}{5}},\quad t > 0,\quad x\in \Bbb{R}.
\label{eq:uxt}
\end{equation}
We can then form the following estimate
\begin{equation}
g_u(u)=(q+1)u^q\leq \tilde{C}(1+t)^{q(1-\frac{q}{5})},\quad t > 0.
\label{eq:guu}
\end{equation}
Note that the exponent can never become larger than $1.25$, which occurs at $q=2.5$. For the second derivative, we have
\begin{equation}
\frac{dg_u}{du}(u)=q(q+1)u^{q-1} \leq\tilde{\tilde{C}}(1+t)^{(q-1)(1-\frac{q}{5})},\quad t>0.
\label{eq:gu}
\end{equation}
If we instead consider a point $s$, close to $u$ we have
\begin{align}
\frac{dg_u}{du}(s)&=q(q+1)s^{q-1}=q(q+1)(u+(s-u))^{q-1}=
q(q+1)\sum_{p=0}^{q-1}\left(\begin{array}{c}q-1\\p\end{array}\right)u^{q-1-p}(s-u)^p\\
&\leq \frac{q!}{\lfloor\frac{q-1}{2} \rfloor!^2}\left\{
\begin{array}{ll}
C(1+t)^{(q-1)(1-\frac{q}{5})}\sum_{p=0}^{q-1}|u-s|^p, & q\leq 5,\\
C(1+t)^{(1-\frac{q}{5})}\sum_{p=0}^{q-1}|u-s|^p, & q>5.
\end{array}\right.
,\quad t>0,
\end{align}
which allows the following estimate
\begin{equation}
|g_u(u)-g_u(s)|\leq C_q(1+t)^{\tilde{q}}\sum_{p=1}^q |u-s|^p,
\label{eq:gugs}
\end{equation}
where $\tilde{q}=\max((q-1)(1-\frac{q}{5}),(1-\frac{q}{5}))$.
%Looking ahead, we will also need to consider integral expressions involving $g_u(u)$ or higher derivatives of $g(u)$ leading to terms of the form
%\[\int_{0}^tf(\tau)(1+\tau)^pd\tau,\]
%where $p\neq -1/2$ is typically not an integer. These integrals are non-trivial to integrate and we will use Cauchy-Schwartz inequality to get
%\begin{equation}
%\left|\int_{0}^tf(\tau)(1+\tau)^pd\tau\right|\leq \left(\int_{0}^t|f(\tau)|^2d\tau\right)^{1/2}\left(\int_{0}^t(1+\tau)^{2p}d\tau\right)^{1/2}=\sqrt{\frac{(1+t)^{2p+1}-1}{2p+1}}\left(\int_{0}^t|f(\tau)|^2d\tau\right)^{1/2}.
%\label{eq:csint}
%\end{equation}
\subsection{Term by term estimates for the error contributions}
In this section, we are going to derive an estimate for the discrete approximation error. We first note that $H(t)$ is a sum over number of terms $H_j(t)$ and split the integral in~\eqref{ode:ersol} to get
\[\|E(t)\|_\infty=\left\|\int_{0}^{t}Q^{-1}H(\tau)d\tau\right\|_\infty
\leq \sum_{j}\left\|\int_{0}^{t}Q^{-1}H_j(\tau)d\tau\right\|_\infty.\]
For the first error term, we have $H_1(t)=-\alpha\epsilon_{xxxx}^\prime(t)$. Integration leads to
\begin{align}
E_1(t)&=-\alpha\int_{0}^{t}Q^{-1}\epsilon^{\prime}_{xxxx}(\tau)d\tau
=-\alpha Q^{-1}\left(\epsilon_{xxxx}(t)-\epsilon_{xxxx}(0)\right),
\label{ode:ersol2}
\end{align}
then we can get the following estimate
\begin{align}
\|E_1(t)\|_\infty&\leq
|\alpha|\|Q^{-1}\|_\infty\max_{0\leq\tau\leq t}\|\epsilon_{xxxx}(\tau)\|_\infty.
\label{eq:e1t}
\end{align}
The second error term that we consider is generated by $H_2(t)=G_u(U)\epsilon_x(t)$, leading to
\begin{align}
\|E_2(t)\|_\infty&=\left\|\int_{0}^{t}Q^{-1}G_u(U)\epsilon_x(\tau)d\tau\right\|_\infty\nonumber\\
&\leq \|Q^{-1}\|_\infty \tilde{C}\max_{0\leq\tau\leq t}\|\epsilon_{x}(\tau)\|_\infty
\int_{0}^{t}(1+\tau)^{q(1-\frac{q}{5})}d\tau\nonumber\\
&=\|Q^{-1}\|_\infty \tilde{C}\max_{0\leq\tau\leq t}\|\epsilon_{x}(\tau)\|_\infty \left(\frac{(1+t)^{q(1-\frac{q}{5})+1}-1}{q(1-\frac{q}{5})+1}\right)
\label{eq:e2t}
\end{align}
where we used~\eqref{eq:guu} for the term involving $G_u(U)$.

The final part focuses on the non-linear term and is the most complicated. We start by rewriting the generating term
\begin{align}
H_3(t)&=G_u(U)\tilde{\Psi}_xU-G_u(S)\tilde{\Psi}_xS+(G_u(U)-G_u(S))F_x(t)\nonumber\\
&=G_u(U)\tilde{\Psi}_x(U-S)+(G_u(U)-G_u(S))(F_x(t)+\tilde{\Psi}_x(U+(S-U))).
\label{eq:nonlin}
\end{align}
If we rewrite equation~\eqref{eq:bondim} in matrix form,
the function $F_x(t)$ can be expressed as
\begin{equation}
F_x(t)=B_xF(t) = \left(\Psi_{x,b}-\Psi_{x,f}B_f^{-1}B_b\quad \Psi_{x,f}B_f^{-1} \right)
\left(\begin{array}{c}F_1(t)\\F_2(t)\end{array}\right).
\end{equation}
Using this, we can provide the first estimate for contribution to the error from the term $H_3(t)$ given by~\eqref{eq:nonlin}
\begin{align}
\|E_3(t)\|_\infty &=
\left\|\int_{0}^{t}Q^{-1}\left(G_u(U)\tilde{\Psi}_x(U-S)+(G_u(U)-G_u(S))(B_xF(\tau)+\tilde{\Psi}_x(U+(S-U)))\right)d\tau\right\|_\infty\nonumber\\
&\leq \|Q^{-1}\|\|\tilde{\Psi}_x\|\int_{0}^{t}\|G_u(U)\| \|U-S\|d\tau\nonumber\\
&+\|Q^{-1}\|\|B_x\|\int_{0}^{t}\|G_u(U)-G_u(S)\|\|F(\tau)\| d\tau\nonumber\\
&+\|Q^{-1}\|\|\tilde{\Psi}_x\|\int_{0}^{t}\|G_u(U)-G_u(S)\|\|U+(S-U)\|d\tau.
\label{ode:ersol3}
\end{align}
Using equations (\ref{eq:uxt}), (\ref{eq:gu}), and (\ref{eq:gugs}) together with the inequality $|U+(S-U)| \leq |U|+|S-U|$ yields
\begin{align}
\|E_3(t)\|_\infty &\leq
\|Q^{-1}\|\|\tilde{\Psi}_x\| \int_{0}^{t}\tilde{C}(1+\tau)^{q(1-\frac{q}{5})} \|U-S\|d\tau\nonumber\\
&+\|Q^{-1}\|\|B_x\|\int_{0}^{t}C_q(1+\tau)^{\tilde{q}}\|F(\tau)\|\sum_{p=1}^{q}\|U-S\|^p d\tau\nonumber\\
&+\|Q^{-1}\|\|\tilde{\Psi}_x\| \int_{0}^{t}C_q(1+\tau)^{\tilde{q}}\sum_{p=1}^{q}\|U-S\|^p\left(C(1+\tau)^{(1-\frac{q}{5})}+\|U-S\|\right)d\tau.
\label{eq:e3t}
\end{align}
%Note that the ${2q(1-\frac{q}{5})+1}$ can never become larger than 3.5, which occurs at $q = 2.5$. Then, we can use the bound $\sqrt{\frac{(1+t)^{2q(1-\frac{q}{5})+1}-1}{2q(1-\frac{q}{5})+1}}\leq C(1+t)^2$. Also,
%we know that all values of $\lambda_j$ are negative. It help us to have the
%$\max_{j}\frac{1-\exp^{-2\beta\lambda_jt}}{2\beta\lambda_j} \leq t$. Therefor the equation~(\ref{eq:termerror}) together equations (\ref{ode:ersol1}) and (\ref{ode:ersol2}) can be rewritten
%where $C_{i,\kappa}$ for $i=1,2,3$ is constant depend on $\kappa_{\infty}(V)$.
Because this term contains the error in the right hand side, it is the most difficult one to include in the error estimate. We simplify it as far as possible. First we write
\begin{align}
\|E_3(t)\|_\infty &\leq \|Q^{-1}\|\max(\|\tilde{\Psi}_x\|,\|B_x\|)
\int_{0}^{t} \sum_{p=1}^{q+1}b_p(\tau)\|E(\tau)\|^p\,\mathrm{d}\tau,
\end{align}
where
\begin{align}
b_1(\tau) &=\tilde{C}(1+\tau)^{q(1-\frac{q}{5})} +C_q(1+\tau)^{\tilde{q}}\left(\|F(\tau)\| + C(1+\tau)^{(1-\frac{q}{5})}\right)\\
b_p(\tau) &=C_q(1+\tau)^{\tilde{q}}\left(\|F(\tau)\|+ C(1+\tau)^{(1-\frac{q}{5})} + 1\right)\\
b_{q+1}(\tau) &=C_q(1+\tau)^{\tilde{q}}
\label{eq:bp}
\end{align}
Then we make the observation that either the error is small and $\|E(\tau)\|\geq\|E(\tau)\|^p$ for $p>1$ or the error is larger than one, in which case $\|E(\tau)\|^{q+1}\geq\|E(\tau)\|^p$ for $p\leq q$. We let $\nu=1$ or $\nu=q+1$, sum up the coefficients and pick the highest power of $(1+\tau)$ to get
\begin{align}
\|E_3(t)\|_\infty &\leq \tilde{C}_q\|Q^{-1}\|\max(\|\tilde{\Psi}_x\|,\|B_x\|)
\int_{0}^{t}(1+\tau)^{\tilde{q}(1-\frac{q}{5})}\|E(\tau)\|^\nu\,\mathrm{d}\tau,
\label{eq:e3tfinal}
\end{align}
where
\begin{equation}
\tilde{C}_q=(q+1)(2+\max_{0\leq\tau\leq t}\|F(\tau)\|)\max(\tilde{C},C_q,C_qC).
\end{equation}
Table~\ref{tab:q} shows the different powers that are involved as a function of $q$. The final power in the estimate grows with $q$.
\begin{table}[!htb]
\centering
\caption{Values of the different powers involved in the error estimates.}
\label{tab:q}
\begin{tabular}{|r|rrrr|}\hline
$q$ &$(1-q/5)$ & $q(1-q/5)$ & $\tilde{q}$ & $\tilde{q}(1-q/5)$ \\\hline
1 & 0.8 & 0.8 & 0.8 & 0.64\\
2 & 0.6 & 1.2 & 0.6 & 0.36\\
3 & 0.4 & 1.2 & 0.8 & 0.32\\
4 & 0.2 & 0.8 & 0.6 & 0.12\\
5 & 0 & 0 & 0 & 0\\
6 & -0.2 & -1.2 & -0.2 & 0.04\\
7 & -0.4 & -2.8 & -0.4 & 0.16\\
8 & -0.6 & -4.8 & -0.6 & 0.36\\
9 & -0.8 & -7.2 & -0.8 & 0.64\\
10 & -1 & -10 & -1 & 1 \\\hline
\end{tabular}
\end{table}

\subsection{Global error bound}
% Check what the first and second terms are
By combining the error terms~\eqref{eq:e1t},~\eqref{eq:e2t} and~\eqref{eq:e3tfinal}, we get the following relation for the error due to the spatial discretization
%In this section, we assumed that the ODE solver is accurate enough and can neglect ODE solver's error in our estimations.
%By collocating the error terms~\eqref{eq:e1t}, \eqref{eq:e2t} and \eqref{eq:e3tfinal}, we can get the following inequality
\begin{align}
\|E(t)\|_{\infty}&\leq C_1\max_{0\leq\tau\leq t}\|\epsilon_{xxxx}(\tau)\|_\infty
+C_2(t)\max_{0\leq\tau\leq t}\|\epsilon_{x}(\tau)\|_\infty \nonumber\\
&+C_3\int_{0}^{t}(1+\tau)^{\tilde{q}(1-\frac{q}{5})}\|E(\tau)\|^\nu\,\mathrm{d}\tau,
\label{eq:termboundF}
\end{align}
where
\begin{align}
C_1& = |\alpha|\|Q^{-1}\|_{\infty},\\
C_2(t)&=\tilde{C}\|Q^{-1}\|_{\infty}\left(\frac{(1+t)^{q(1-\frac{q}{5})+1}-1}{q(1-\frac{q}{5})+1}\right),\\
C_3 &=\tilde{C}_q\|Q^{-1}\|\max(\|\tilde{\Psi}_x\|,\|B_x\|).
\label{eq:consts}
\end{align}
%\begin{align}
%\|E(t)\|_{\infty}&\leq C_{1}\max_{0\leq\tau\leq t}\|\epsilon_{xxxx}(\tau)\|_\infty
%+C_2(t)\max_{0\leq\tau\leq t}\|\epsilon_{x}(\tau)\|_\infty \nonumber\\
%&+C_{3}\int_{0}^{t}\left(k(\tau)\|E(\tau)\|+\ldots+(1+\tau)^{\tilde{q}}\|E(\tau)\|^{q+1}\right)\,\mathrm{d}\tau,
%\label{eq:termboundF}
%\end{align}
%where $C_1=|\alpha|\|Q\|_{\infty}$ and $C_2(t)=\tilde{C}\|Q\|_{\infty}\left(\frac{(1+t)^{q(1-\frac{q}{5})+1}-1}{q(1-\frac{q}{5})+1}\right)$.
To convert this into an error estimate, we need the following Gronwall inequality~\cite{drago}:
%Now, we need the Gronwall lemma to estimate the global error bound on the discretization points.
\begin{lemma}[Gronwall inequality]
Let the functions $E(t)$, $a(t)$ and $k(t)\geq 0$ be continuous functions defined for $t\in[0,b]$. We assume that for $t\in[0,b]$ we have the inequality
\begin{equation}
E(t)\leq a(t)+\int_{0}^{t}k(\tau)E^n(\tau)\, \mathrm{d}\tau.
\end{equation}
Then for the case $n=1$ it holds
\begin{equation}
E(t)\leq a(t)+\int_{0}^{t}k(\tau)a(\tau) e^{\int_{\tau}^{t}k(u)du}\, \mathrm{d}\tau.
\end{equation}
If the function $a(t)$ is also non-decreasing, then
\begin{equation}
E(t)\leq a(t) e^{\int_{0}^{t}k(\tau)\, \mathrm{d}\tau}.
\end{equation}
For the case $n\geq 2$
\begin{equation}
E(t)\leq a(t)\left[1-(n-1)\int_{0}^{t}k(\tau)a^{n-1}(\tau)\, \mathrm{d}\tau\right]^{\frac{1}{n-1}}, \quad t\leq \beta_n,
\label{eq:nongra}
\end{equation}
where $\beta_n=\sup\{t:(n-1)\int_{0}^{t}k(\tau)a^{n-1}(\tau)\, \mathrm{d}\tau<1\}$.
\label{lem:gran}
\end{lemma}
In our case, it can easily be verified that the function
\begin{equation}
a(t)=C_1\max_{0\leq\tau\leq t}\|\epsilon_{xxxx}(\tau)\|_\infty
+C_2(t)\max_{0\leq\tau\leq t}\|\epsilon_{x}(\tau)\|_\infty
\end{equation}
is non-decreasing in time. For the case of small errors, $n=\nu=1$, and
\begin{equation}
\int_{0}^{t}k(\tau)\, \mathrm{d}\tau=C_3\int_{0}^{t}(1+\tau)^{\tilde{q}(1-\frac{q}{5})}\, \mathrm{d}\tau = C_3\frac{(1+t)^{\tilde{q}(1-\frac{q}{5})+1}-1}{\tilde{q}(1-\frac{q}{5})+1},
%=\tilde{C}_3\left((1+t)^{\tilde{q}(1-\frac{q}{5})+1}-1\right)
\end{equation}
%
%If we consider cases where the error $\|E(t)\|_{\infty}$ is smaller than one and ignoring the error items with the powers higher than one in the right hand side of equation~\eqref{eq:termboundF}, applying
the Gronwall inequality leads to
\begin{align}
\|E(t)\|_{\infty}\leq
\left[C_{1}\max_{0\leq\tau\leq t}\|\epsilon_{xxxx}(\tau)\|_\infty
+C_2(t)\max_{0\leq\tau\leq t}\|\epsilon_{x}(\tau)\|_\infty
\right] e^{C_3\frac{(1+t)^{\tilde{q}(1-\frac{q}{5})+1}-1}{\tilde{q}(1-\frac{q}{5})+1}}.
\label{eq:diserror}
\end{align}
For the case of errors larger than one, we do not carry out the full derivation, but note that the limit on the time interval becomes less severe when $a(t)$ is small enough, which is the case when the spatial resolution is high enough.
%\subsection{Estimates for the RBF interpolation error}
%%
%For any $k\in{\Bbb{N}}$, Sobolev space on $\Bbb{R}^d$ is defined
%by
%\begin{equation}
%W_2^k(\Bbb{R}^d)=\{f\in L_2(\Bbb{R}^d):
%\hat{f}(\cdot)(1+\|\cdot\|_2^2)^{\frac{k}{2}}\in L_2(\Bbb{R}^d)\},
%\end{equation}
%where $\hat{f}$ is Fourier transform of $f$.
%
%Let $\hat{\phi}$ be the Fourier transform of the inverse multiquadric radial basis function $\phi$. The native Hilbert for the case of inverse multiquadric kernel is defined as
%\begin{equation}
%\mathcal{N}_{\phi}=\{f\in C(\Bbb{R}^d)\cap L_2(\Bbb{R}^d):
%\frac{\hat{f}}{\sqrt{\hat{\phi}}}\in L_2(\Bbb{R}^d)\}.
%\end{equation}
%From \cite{rigzwi}, we know for all non-negative $k$ and a bounded domain $\Omega$, native space $\mathcal{N}_{\phi}(\Omega)$ is contained in Sobolev space $W_2^k(\Omega)$. i.e., $\mathcal{N}_{\phi}(\Omega)\subset W_2^k(\Omega)$. Furthermore
%there are constants $c_1$ and $c_2$ such that
%\begin{equation}
%\|f\|_{W_2^k(\Omega)}\leq \max\left\{\frac{2^k}{c_1},c_2^kk^{2k}\right\}^{\frac{1}{2}}\|f\|_{\mathcal{N}_{\phi}},~ \textit{for all}~ k\in\Bbb{N}.
%\end{equation}

% It remains to find interpolation error estimates. We have x_j in 1-D, but we need to keep also 2-D as a general case.
It remains to insert estimates for the derivatives of the RBF interpolation errors. These are typically measured in terms of the fill distance $h$, which in one space dimension with uniform nodes becomes $h=\frac{1}{2}|x_{j+1}-x_j|$, and more generally for a domain $\Omega$ in $d$ dimensions and a node set $X$ is defined as
%
%We define the density of a given set $X_{\Omega}$ consisting of
%pairwise distinct points in $\Omega$ by fill distance
\begin{equation}
h=\sup_{x\in\Omega} \min_{x_j\in
X}\|x-x_j\|.
\end{equation}
Exponential estimates for inverse multiquadric interpolants are given in~\cite{rigzwi} provided that $\Omega$ is a bounded domain with Lipschitz boundary that satisfies an interior cone condition.
% Swap to ep instead of what we have now.
%Let $\Omega\subset\Bbb{R}^d$ is bounded with Lipschitz boundary and satisfies an interior cone condition. Also, we assume that $\mathcal{L}$ be the
%linear differential operator of order $m$, and $u$ is the exact solution and $z$ is the RBF interpolant on the domain $\Omega$. Then we can use the following estimates from~\cite{rigzwi} for the inverse multiquadrics function interpolation errors and their derivatives
\begin{equation}
\|\epsilon_\mathcal{L}\|_\infty\leq
c e^{-{\gamma/\sqrt{h}}}\|u\|_{\mathcal{N}(\Omega)},
\label{eq:intererror}
\end{equation}
where the constant $\gamma$ depends on the order of the differential operator $\mathcal{L}$, the dimension $d$, and the shape parameter of the RBF $\ep{}$, and $c$ is a positive constant. The norm in the right hand side denoted by $\|.\|_{\mathcal{N}}(.)$ is the native space norm. For more details about this norm see~\cite{Fass,rigzwi}.

Now, by inserting the interpolation error estimate~\eqref{eq:intererror} into the global error estimate~\eqref{eq:diserror}
we get
\begin{align}
\|E(t)\|_{\infty}\leq C(t) e^{-\frac{\gamma}{\sqrt{h}}}e^{C_3\frac{(1+t)^{\tilde{q}(1-\frac{q}{5})+1}-1}{\tilde{q}(1-\frac{q}{5})+1}} \max_{0\leq\tau\leq t}\|u(\tau)\|_{\mathcal{N}(\Omega)}.
\label{eq:gdiserror}
\end{align}
where $C(t)=c(C_1+C_2(t))$. The final estimate tells us that the spatial RBF discretization does allow for exponential convergence in $h$, but we should expect the error to grow in time. For any finite interval $t\in[0,T]$ the growth in time is however bounded.

%Remark: we know that term $C_{1}\max_{0\leq\tau\leq t}\|\epsilon_{xxxx}(\tau)\|_\infty
%+C_2(t)\max_{0\leq\tau\leq t}\|\epsilon_{x}(\tau)\|_\infty$ in equation~\eqref{eq:termboundF} is a function of interpolation error and of course for the enough small value of the $h$ can expect to be the so small value. Therefore, for all values of $t$ appropriate to the not so small value of $\beta_n$, from the Gronwall inequality for $n\geq 2$, we can get
%\begin{equation}
%\|E(t)\|_{\infty}\leq C_{1}\max_{0\leq\tau\leq t}\|\epsilon_{xxxx}(\tau)\|_\infty
%+C_2(t)\max_{0\leq\tau\leq t}\|\epsilon_{x}(\tau)\|_\infty.
%\end{equation}
%It means that if we consider the error items with the powers higher than one in right hand side of equation~\eqref{eq:termboundF}, we will have
%\begin{align}
%\|E(t)\|_{\infty}\leq c C(t) e^{-\frac{\gamma}{\sqrt{h}}} \max_{0\leq\tau\leq t}\|u(\tau)\|_{\mathcal{N}(\Omega)}.
%\end{align}

\subsection{Numerical investigation of the matrix norms in the estimates}
\begin{figure}[h!]
\centering
\includegraphics[width=0.4\textwidth]{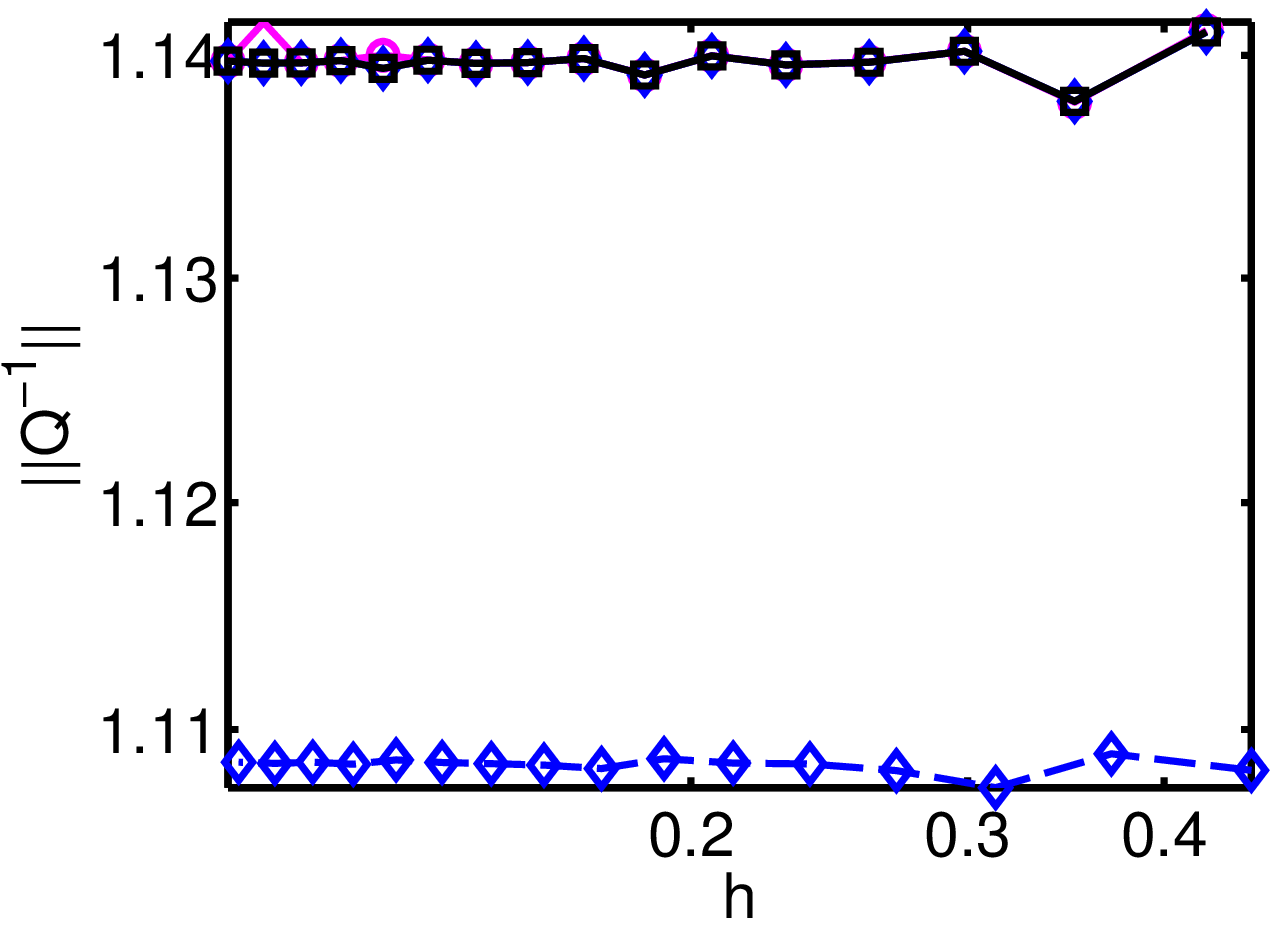}\\
\includegraphics[width=0.4\textwidth]{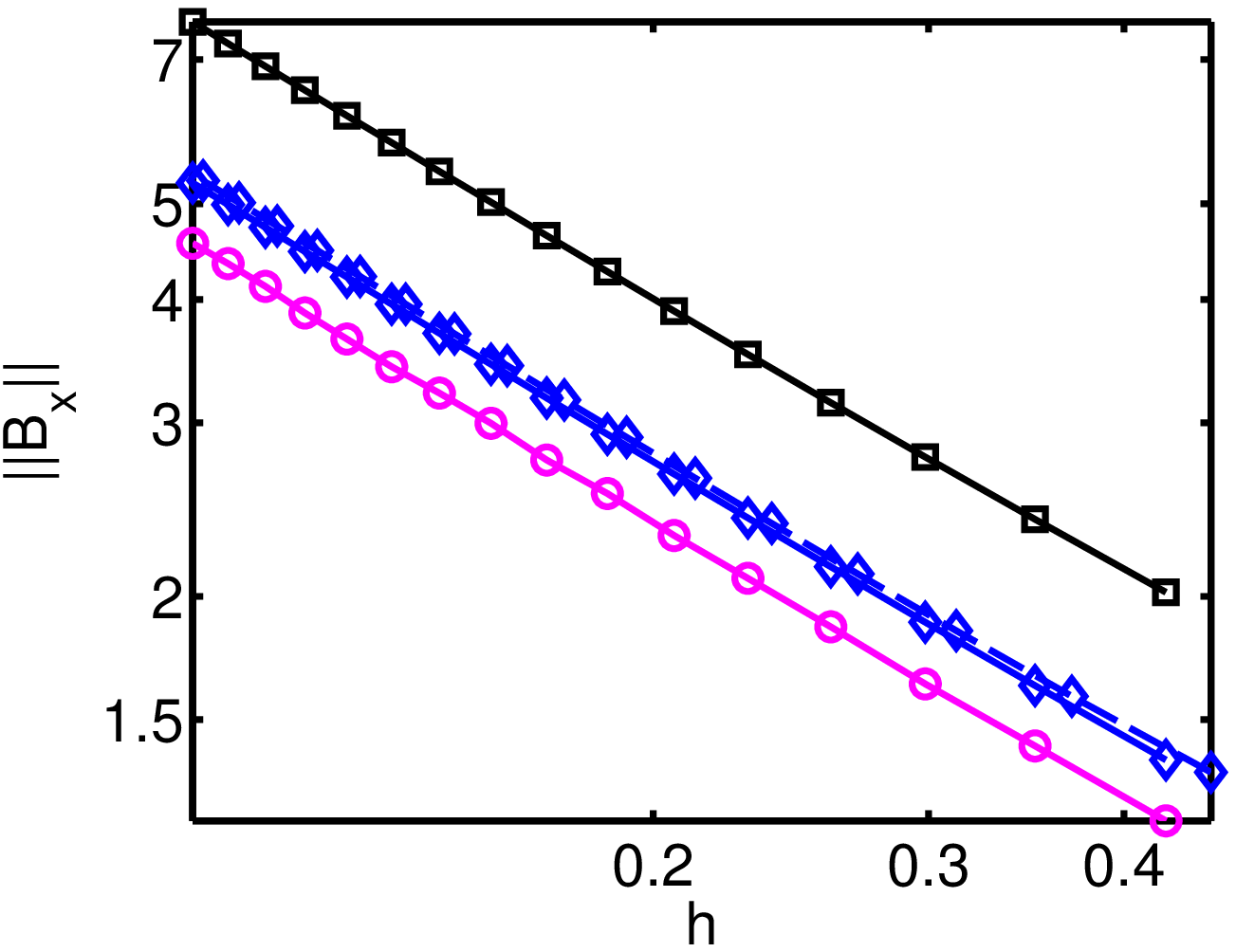}
\includegraphics[width=0.4\textwidth]{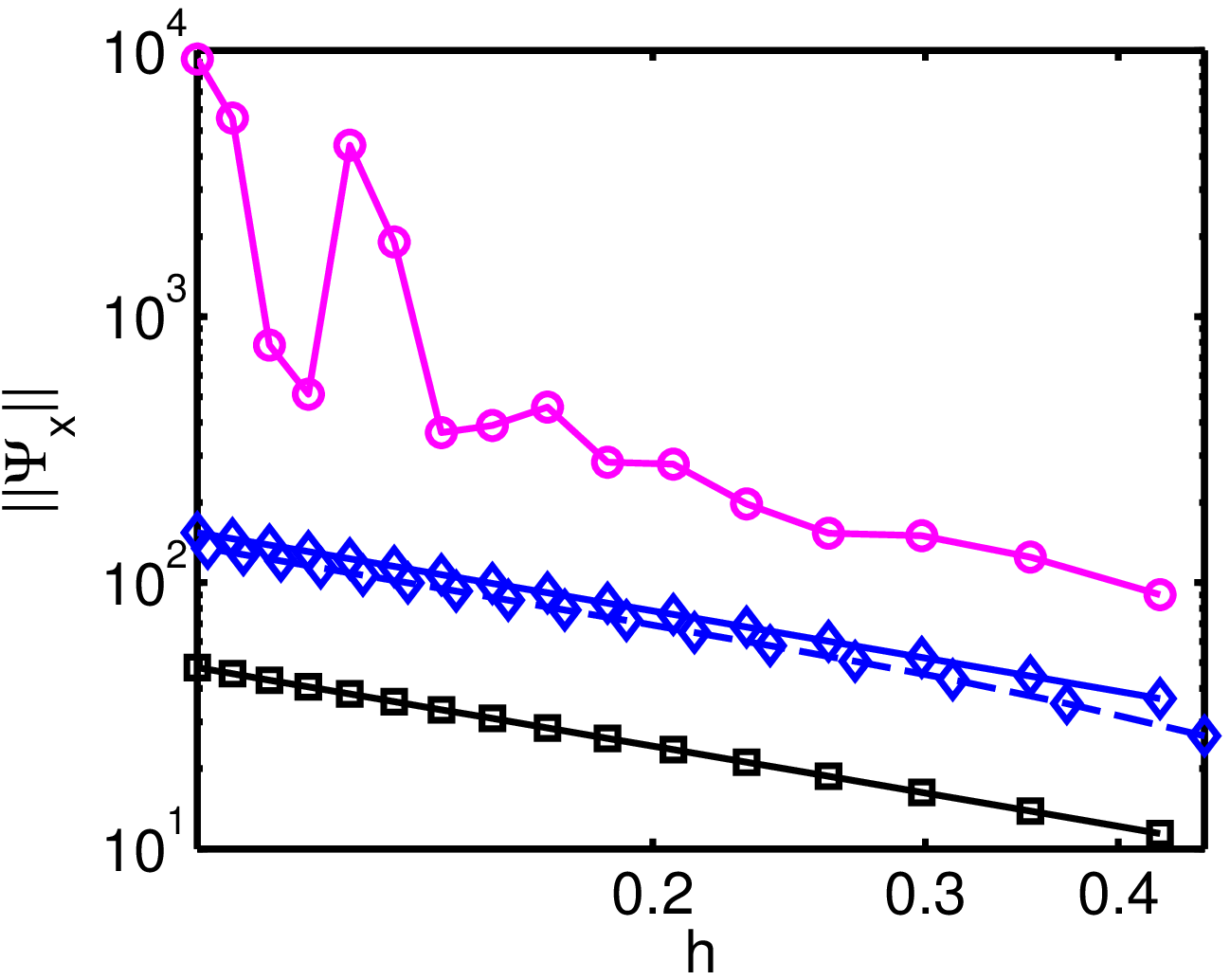}
\caption{The norm $\|Q^{-1}\|_{\infty}$ (top), the norm $\|B_x\|_{\infty}$ (bottom left), and the norm $\|\tilde{\Psi}_x\|_{\infty}$ (bottom right) as a function of $h$ for $\ep{}h=0.4$ ($\circ$), $\ep{}h=0.5$ ($\diamond$), $\ep{}h=1$ ($\Box$) and $L=10$ (solid lines) and $L=5$ (dashed lines). }
\label{Qpar}
\end{figure}

There are three different matrix norms that appear in the error estimates. We do not have theoretical bounds for these, and therefore we perform a numerical investigation of their behaviors. Based on previous experience of RBF approximation, we have selected the following representation of the method parameters, the fill distance $h$, the relative shape parameter value $\ep{}h$, and the (half) domain size $L$. In Figure~\ref{Qpar}, the norms are plotted as a function of $h$ for different combinations of the parameters. The chosen values are such that they are reasonable for approximation. By choosing extreme values, different behaviors can be observed. We see that the norm $\|Q^{-1}\|$ does not change much with $h$ or $\ep{}h$, but slightly with $L$. Both of the norms $\|B_x\|$ and $\|\tilde{\Psi}_x\|$ grow as $h^{-1}$, and we can observe a little bit of instability in the value for small $\ep{}$. This rate is what would be expected for a first derivative such as is represented by these matrices. Changing $L$ has a very small effect also in this case.

%These norms that appear in all the terms are not scale invariant. Therefore,
%The numerical investigation on figure~\ref{Qpar} shows how $\|Q\|_{\infty}\|\tilde{\Psi}_x\|_{\infty}$ and $\|Q\|_{\infty}\|B\|_{\infty}$ varies with the RBF approximation parameters when the domain $\Omega$ is unit. All cases in the figure~\ref{Qpar} are bounded to one, although the maximum value is belong to $\|Q\|_{\infty}\|\tilde{\Psi}_x\|_{\infty}$. However the differences are not large in the range of $\ep{}$-values when explored without running into ill-conditioning.

%\comment{need explaining the $\Omega=[-L,L]$ and scaling with length of domain.}
% Question, how does mapping of the interval affect the non-linear term? y=x/L du/dy=dx/dy*du/dx = L du/dx alpha/L^4 och g_u/L
%

\section{\textsc{MATLAB} Implementation}
\label{sec:matlab}
In this section, sample MATLAB implementations of the fictitious point and resampling RBF methods for the one-dimensional Rosenau equation~(\ref{eq:Rosenau1D})--(\ref{eq:1Dic}) are presented and discussed. We use an example with a known solution. For $g(u)=10u^3-12u^5-\frac{3}{2}u$ and $\alpha(x,t)=0.5$ it holds that $u(x,t)=\textnormal{sech}(x-t)$ is a solution~\cite{park}. For both methods, equally spaced nodes are used, and the spatial domain is $[-L,L]$. 
%The Rosenau equation (\ref{eq:Rosenau1D})--(\ref{eq:1Dic}) with $g(u)=10u^3-12u^5-\frac{3}{2}u$ and $\alpha(x,t)=0.5$ is considered
%with the known exact solution $u(x,t)=\textnormal{sech}(x-t)$.
%Equally-spaced nodes is used to discretize the interval $[-L,L]$.

\subsection{Implementation of the fictitious point method}
\begin{figure}[h!]
\centerline{
\includegraphics[width=0.55\textwidth]{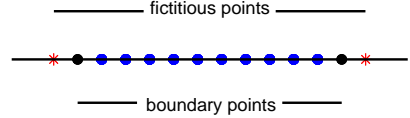}}
\caption{An example of a node distribution for the  fictitious point method in the one-dimensional case.}
\label{fig:node}
\end{figure}

Following the idea of the fictitious point method in subsection~\ref{sec:fp}, we complement the interior and boundary RBF nodes with two (the number of boundary nodes) fictitious points outside the computational domain, see Figure~\ref{fig:node} for an illustration.
We generate the differentiation matrices using the modified basis functions according to Equation~\eqref{eq:moddiff}. 
%and then spilt the differentiation matrices into separate compartment corresponding to the interior, boundary and fictitious points. 
Collocating the Rosenau equation by applying the modified differentiation matrices leads to ODE system~\eqref{eq:fpmol}, which we here solve by using the built-in MATLAB ODE solver \textbf{ode15s}. The two functions below constitute a complete MATLAB implementation of the problem.
{\small
\begin{verbatim}
function [S,x,t]=fictitious(N,L,T)
% N : The number of node points
% L : [-L,L] is the domain
% T : Final time
% Exact solution and derivatives for test case  
  u = @(x,t) sech(x-t); 
  ux = @(x,t) sech(t-x).*tanh(t-x); 
  ut = @(x,t) -sech(t-x).*tanh(t-x);
  uxt = @(x,t) sech(t - x).^3 - sech(t - x).*tanh(t - x).^2;
% Generate nodes x. Map x to [-L,L] such that x(2) = -L, x(N-1) = L
% and x(1), x(N) are the left and right fictitious points respectively.
  x = linspace(-L,L,N); 
  linmap = @(x,x1,x2,y1,y2) (y2-y1)*(x-x1)/(x2-x1) + y1;
  x = linmap(x,x(2),x(N-1),-L,L); x = x(:); 
% Differentiation matrices for inverse multiquadric RBF  
  phi = @(ep,r) 1./sqrt(1+(ep*r).^2); 
  ep = 0.08/min(diff(x)); % Shape parameter
  dx = bsxfun(@minus,x,x.'); A = phi(ep,dx);
  Dx = (-ep^2*dx.*A.^3)/A; % 1st Derivative matrix
  D4x = (3*ep^4*(3-24*(ep*dx).^2+8*(ep*dx).^4).*A.^9)/A; % 4th Derivative
  Bf = Dx([2 N-1],[1 N]); Bd = Dx([2 N-1],3:N-2); Bb = Dx([2 N-1],[2,N-1]);
% Modify differentiation matrices
  Dxd = Dx(3:N-2,3:N-2) - (Dx(3:N-2,[1 N])/Bf)*Bd;
  Dxb = [Dx(3:N-2,[2 N-1]) - (Dx(3:N-2,[1 N])/Bf)*Bb  Dx(3:N-2,[1 N])/Bf];
  D4xd = D4x(3:N-2,3:N-2) - (D4x(3:N-2,[1 N])/Bf)*Bd;
  D4xb = [D4x(3:N-2,[2 N-1]) - (D4x(3:N-2,[1 N])/Bf)*Bb  D4x(3:N-2,[1 N])/Bf];
% Initial condition
  S0 = u(x(3:N-2),0); opt = odeset('RelTol',1e-10);
% Solve the ODE for the approximate solution S
  fun = @(t,S) odefun(t,S,x,N,Dxd,Dxb,D4xd,D4xb,u,ux,ut,uxt);
  [t,S] = ode15s(fun,[0 T],S0,opt);
% Plot the solution for all times 
  figure(1),clf,plot(x(3:N-2),S)
% Plot the error at the final time
  E = abs(S(end,:)-u(x(3:N-2),T)');
  figure(1),clf,plot(x(3:N-2),E)
  
function Sprime = odefun(t,S,x,N,Dxd,Dxb,D4xd,D4xb,u,ux,ut,uxt)
  Fx = [u(x([2 N-1]),t); ux(x([2 N-1]),t)];
  F4xt = [ut(x([2 N-1]),t); uxt(x([2 N-1]),t)];
  Gu = diag(-1.5 - 60*S.^4 + 30*S.^2);
  Sprime = (eye(N-4) + 0.5*D4xd)\(Gu*Dxd*S + Gu*Dxb*Fx - 0.5*D4xb*F4xt);
\end{verbatim}
}
It can be noted that when we use the modified basis functions, we need to provide the time derivatives of the boundary conditions as well as the boundary conditions themselves. This is not needed with the alternative formulation~\eqref{eq:fict}, but instead the system is stated in DAE form.

\subsection{ Implementation of the resampling RBF method}
For the resampling method, we start directly from the DAE form derived in Section~\ref{sec:res}, where the four boundary conditions enforced at the boundary points constitute the algebraic part. The $N-4$ auxiliary points where the PDE is enforced are uniformly distributed in the interior of the computational domain and do not in general coincide with the RBF node points where the solution is sought. We organize the solution vector as $S=[S_d ~S_b]^T$, where as before, $S_d$ contains solution values at the interior RBF nodes, and $S_b$ contains the two boundary values. Then the DAE discretization scheme can be schematically be displayed  as
\begin{equation*}
\begin{tabular}{lllll}
% First row
\fblockmatrix{0.75in}{0.5in}{$\Psi^R + \frac{1}{2}\Psi^R_{xxxx}$}
\fblockmatrix{0.12in}{1.32in}{$S^\prime$} &
\blockmatrix{0.1in}{1.2in}{$=$} & \fblockmatrix{1in}{0.5in}{$G_u(\tilde{S})\Psi^R_x$}
\fblockmatrix{0.1in}{1.32in}{$S$} & \blockmatrix{0.1in}{1.2in}{$-$}
& \fblockmatrix{0.35in}{0.5in}{$0$}
\vspace*{-.82in}
\\
% Second row
\fblockmatrix{0.75in}{0.25in}{$0$} & &
\fblockmatrix{1in}{0.25in}{$
\begin{matrix}
0&\cdots&0&1&0\\
0&\cdots&0&0&1\\
\end{matrix}
$} & & \fblockmatrix{0.35in}{0.25in}{$F_1(t)
%\begin{matrix}
%f_1(t)\\
%f_2(t)\\
%\end{matrix}
$}
\\
%Third row
\fblockmatrix{0.75in}{0.25in}{$0$} & &
\fblockmatrix{1in}{0.25in}{$
\tilde{B}
$} & & \fblockmatrix{0.35in}{0.25in}{$ F_2(t)$}
\end{tabular},
\end{equation*}
where $\Psi^R$ is an $N-4 \times N$ resampling matrix that provides values at the auxiliary points given values at the node points, $\tilde{S}=\Psi^RS$, $G_u(\tilde{S})$ is an $(N-4) \times (N-4)$ diagonal matrix, $\Psi^R_{x}$ and $\Psi^R_{xxxx}$ are $N-4 \times N$ resampled first and
fourth derivative matrices respectively. The derivative boundary condition matrix $\tilde{B}=(\tilde{B}_d~~\tilde{B}_b)$ is a $2\times N$ matrix, see Equation~\eqref{eq:resamplingbc} for details.
%
%The matrix operator $\Psi^R$ can be also seen as ``Resampling" or ``Projection" matrix operator.
%In the diagram above, $\Psi^R$ is applied to $S$ first resulting in interpolation values at the $N-4$ auxiliary points and then enforce the PDEs at those points.

%However,
%one may prefer to change the order by applying $\Psi^R$ to $S$ first resulting in interpolation values at the
%$N-4$ auxiliary points and then enforce the PDEs at those points.
% In our case, we do not see any difference in accuracy between the two implementations.

Both ODEs and DAEs of index 1 can be solved in MATLAB using \textbf{ode15s}, previously employed for the fictitious point method. One may also use the syntactically similar open source software \textsc{OCTAVE} and there use \textbf{dasspk} as DAE solver. The following two functions show the MATLAB implementation of the resampling RBF method.
{\small
\begin{verbatim}
function [S,x,t]=resampling(N,L,T)
% N : The number of node points
% L : [-L,L] is the domain
% T : Final time
% Exact solution and derivatives for test case  
  u = @(x,t) sech(x-t); 
  ux = @(x,t) sech(t-x).*tanh(t-x);   
% Generate N uniform RBF nodes with boundary pts last
  x = linspace(-L,L,N).'; 
  x = [x(2:end-1); x([1 end])];
% Generate N-4 uniform auxiliary interior points
  xt = linspace(-L,L,N-3).'; xt(end) = []; xt = xt + 0.5*min(diff(xt));
% Differentiation matrices for inverse multiquadric RBF  
  phi = @(ep,r) 1./sqrt(1+(ep*r).^2); 
  ep = 0.08/(x(2)-x(1)); % Shape parameter
  r = bsxfun(@minus,x,x.'); A = phi(ep,r);
% First derivative matrix at x = -L and x = L
  r = bsxfun(@minus,x([N-1 N]),x.');
  Bt = (-ep^2*r.*phi(ep,r).^3)/A;
% Rectangular projection from x to xt
  r = bsxfun(@minus,xt,x.'); R = phi(ep,r);
  PRx = (-ep^2*r.*R.^3)/A; PR = R/A;
  PR4x = (3*ep^4*(3-24*(ep*r).^2+8*(ep*r).^4).*R.^9)/A;
% Initial condition
  S0 = u(x,0); 
  M = [(PR + 0.5*PR4x); zeros(4,N)];
  opt = odeset('mass',M,'masssing','yes','RelTol',1e-9);
  [t,S] = ode15s(@(t,S) daefun(t,S,x,N,PR,PRx,Bt,u,ux),[0 T],S0,opt);
% Plot the solution for all times
  ind = [N-1 2:N-2 N];
  figure(1),clf,plot(x(ind),S(:,ind))
% Plot the error at the final time
  E = abs(S(end,:)-u(x,T)');
  figure(2),clf,plot(x(ind),E(ind))
  
function Sprime = daefun(t,S,x,N,PR,PRx,Bt,u,ux)
  F = [zeros(N-4,1); u(x([N-1 N]),t); ux(x([N-1 N]),t)];
  Gu = diag(- 60*(PR*S).^4 + 30*(PR*S).^2 - 1.5);
  Sprime = [Gu*PRx; [zeros(2,N-2) eye(2,2)]; Bt]*S - F;

\end{verbatim}
}
\section{Numerical results}\label{sec:num}

In this section, we perform numerical experiments to investigate the accuracy and convergence of the proposed schemes. Both one-dimensional and two-dimensional test cases are considered. In all tests, the inverse multiquadric RBF is used. The shape parameter has not been optimized for accuracy. Instead, the parameter range has been chosen such that ill-conditioning is avoided. For the one-dimensional test case, we compare the results with those of a pseudo-spectral resampling method. We have not included the code here, but it can be downloaded from the authors' web pages. 

\subsection{The one-dimensional case}
We consider the same test problem~\cite{park} for the Rosenau equation~\eqref{eq:Rosenau1D} as in the previous section, with the known exact solution $u(x,t)=\textrm{sech}(x-t)$, obtained for $g(u)=10u^3-12u^5-\frac{3}{2}u$ and $\alpha(t)=0.5$. 
The initial and boundary conditions are taken from the exact solution.

%As an example, we consider the Rosenau equation given in
%\cite{park}:
%\begin{align}
%u_t+\frac{1}{2}u_{xxxxt}=g(u)_x,\label{eq:test}
%\end{align}
%where $g(u)=10u^3-12u^5-\frac{3}{2}u$. The exact solution for the
%equation~(\ref{eq:test}) is known to be
%$u(x,t)=\textrm{sech}(x-t)$. We take the value of the exact
%solution at $t=0$ as our initial condition
%\begin{align*}
%u(x,0)=\textrm{sech}(x),
%\end{align*}
%and boundary conditions
%\begin{align*}
%u(-L,t)&=\textrm{sech}(-L-t),&u&(L,t)=\textrm{sech}(L-t),\\
%u_x(-L,t)&=-\textrm{sech}(-L-t)\textrm{tanh}(-L-t),&u_x&(L,t)=-\textrm{sech}(L-t)\textrm{tanh}(L-t).
%\end{align*}

\begin{figure}[!htb]
\centering
\includegraphics[width=0.49\textwidth]{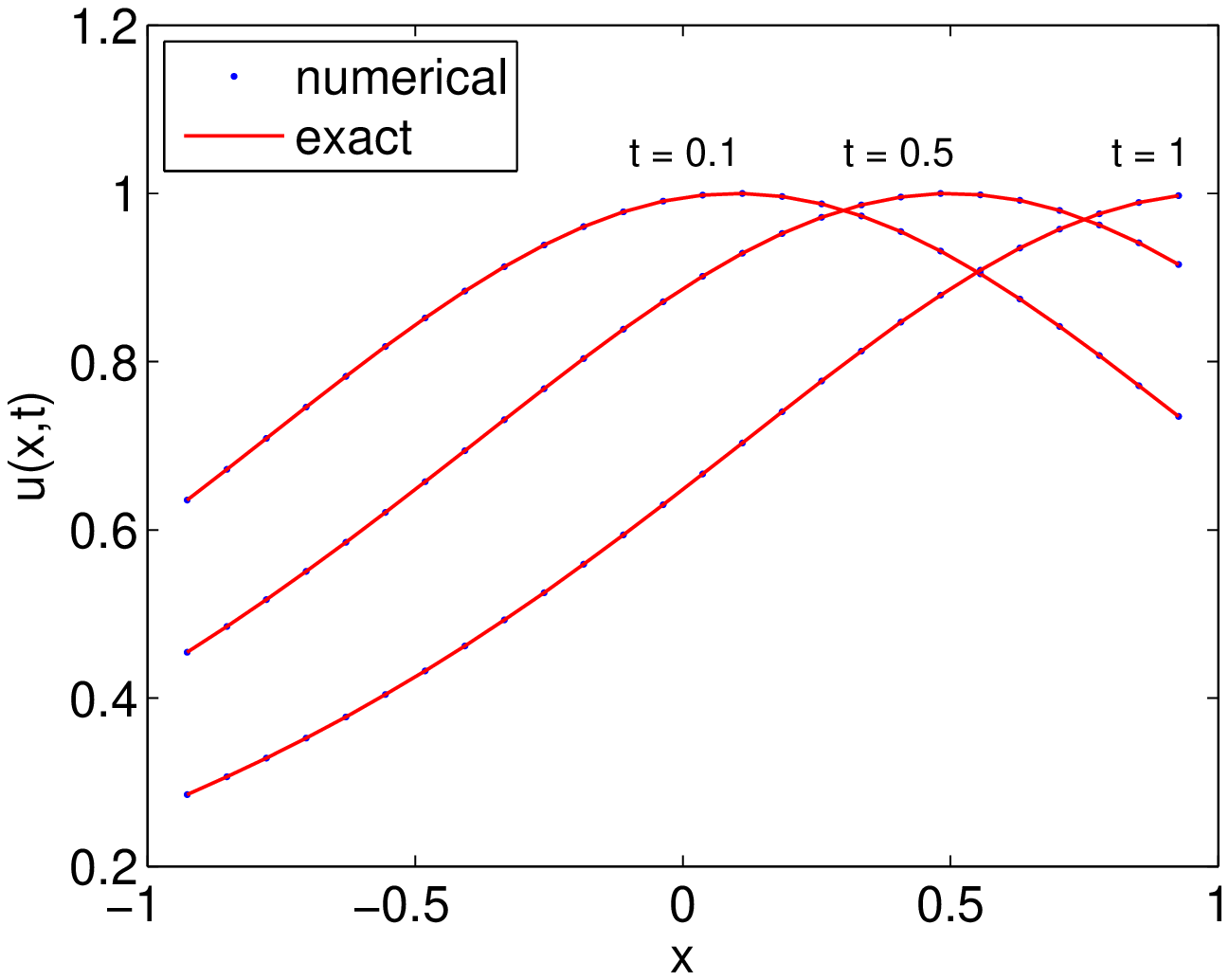}
\includegraphics[width=0.49\textwidth]{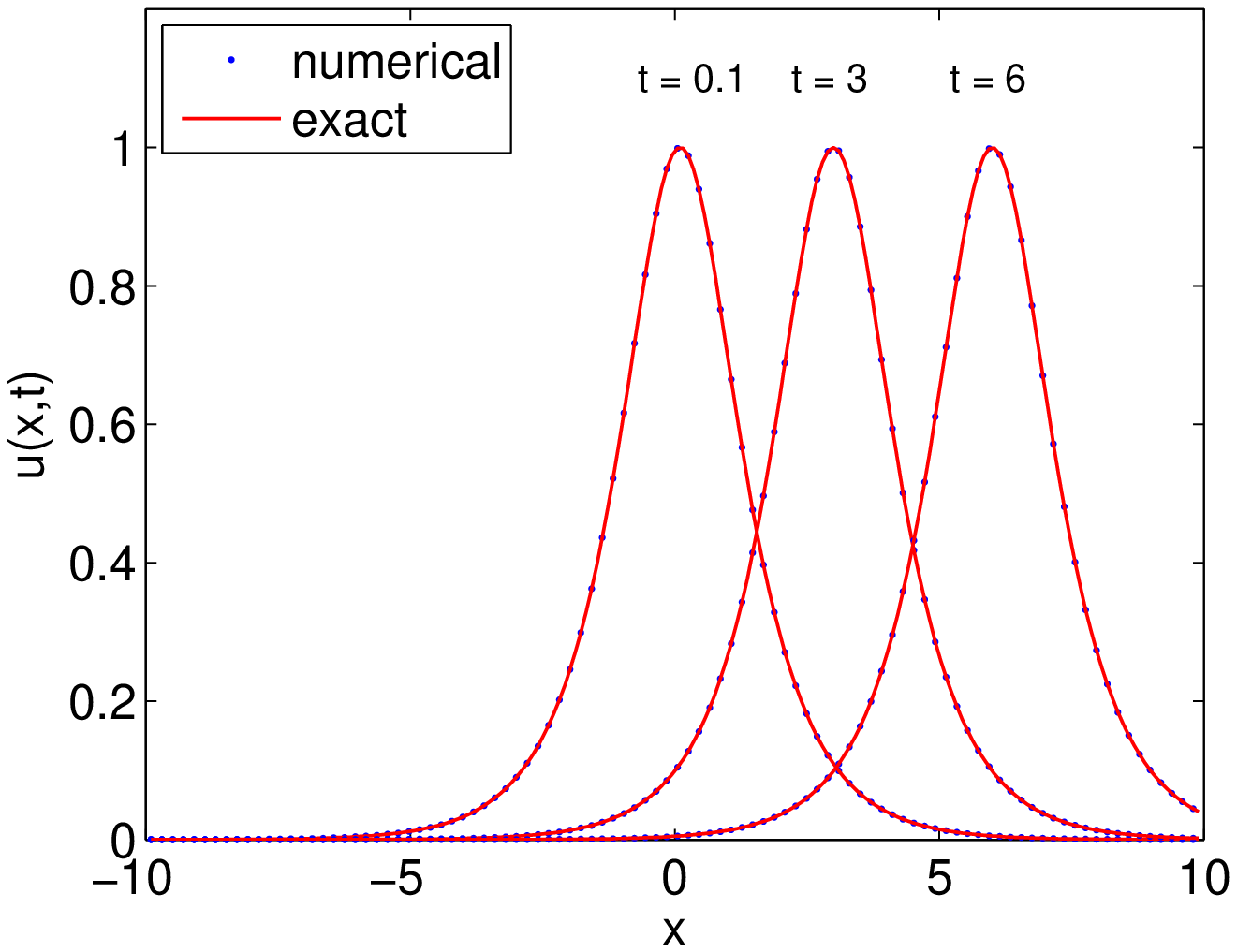}
\caption{Exact and numerical solutions based on the fictitious point method with $L=1$ and $N=30$ (left) and $L=10$ and $N=100$ (right) with $\ep{}=\frac{0.08}{h}$.}
\label{fig:numsol}
\end{figure}

The exact solution over the interval $[-L,\,L]$ for $L=1$ and $L=10$ is shown at different times $t$ together with the numerical solution using the fictitious point method in Figure~\ref{fig:numsol}. The solution is a pulse that travels to the right with time. 

\begin{figure}[h!]
\centering
\includegraphics[width=0.49\textwidth]{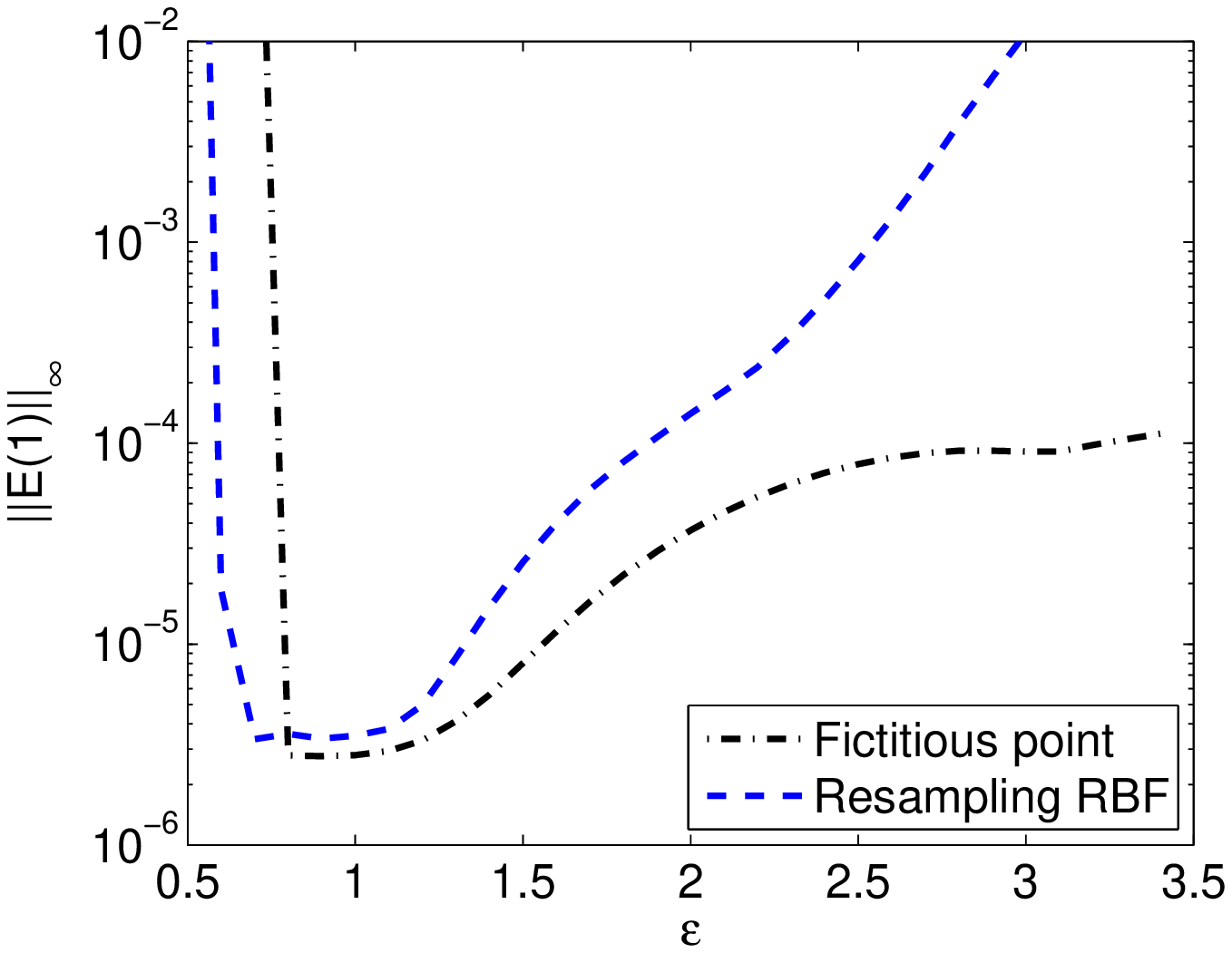}
\includegraphics[width=0.49\textwidth]{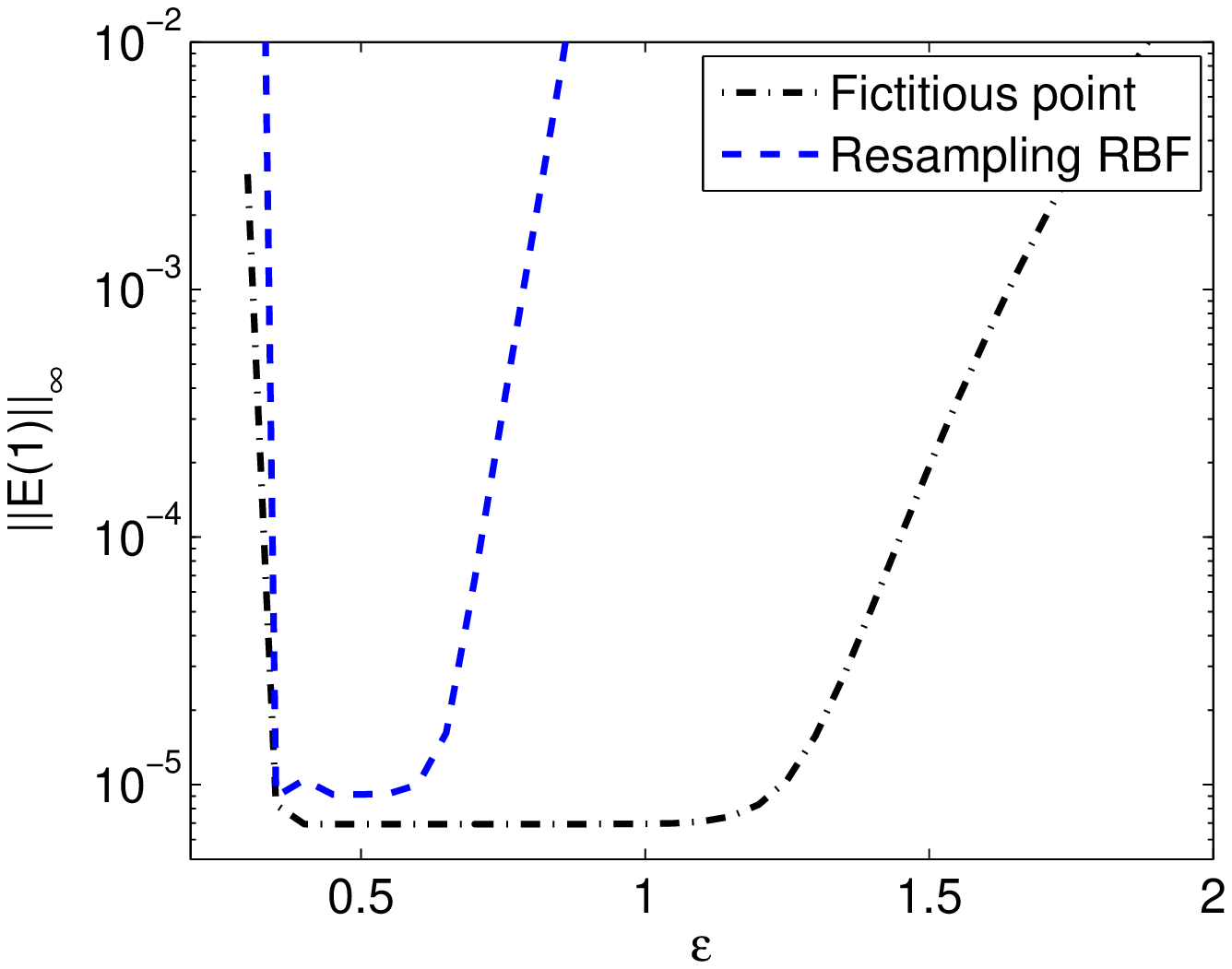}
\caption{The $L_{\infty}$ error at time $t=1$ as a function of the shape parameter $\ep{}$ for $L=1$ and $N=30$ (left) and $L=10$ and $N=100$ (right).}
\label{fig:ep}
\end{figure}

In Figure~\ref{fig:numsol}, a shape parameter $\ep{}=\frac{0.08}{h}=\frac{0.04(N-1)}{L}$ is used. This relation is experimentally determined to ensure stable computations and highest possible solution accuracy.
Figure~\ref{fig:ep} shows how the errors of the two RBF methods vary with $\ep{}$. Using the formula leads to $\ep{}=1.16$ and $\ep{}=0.4$ for the two cases, which is within the best region for each method. It can be noted that the stable range of $\ep{}$ is narrower for the resampling method and that both methods rapidly become ill-conditioned when $\ep{}$ is too small.

\begin{figure}[htb]
\centering
\includegraphics[width=0.49\textwidth]{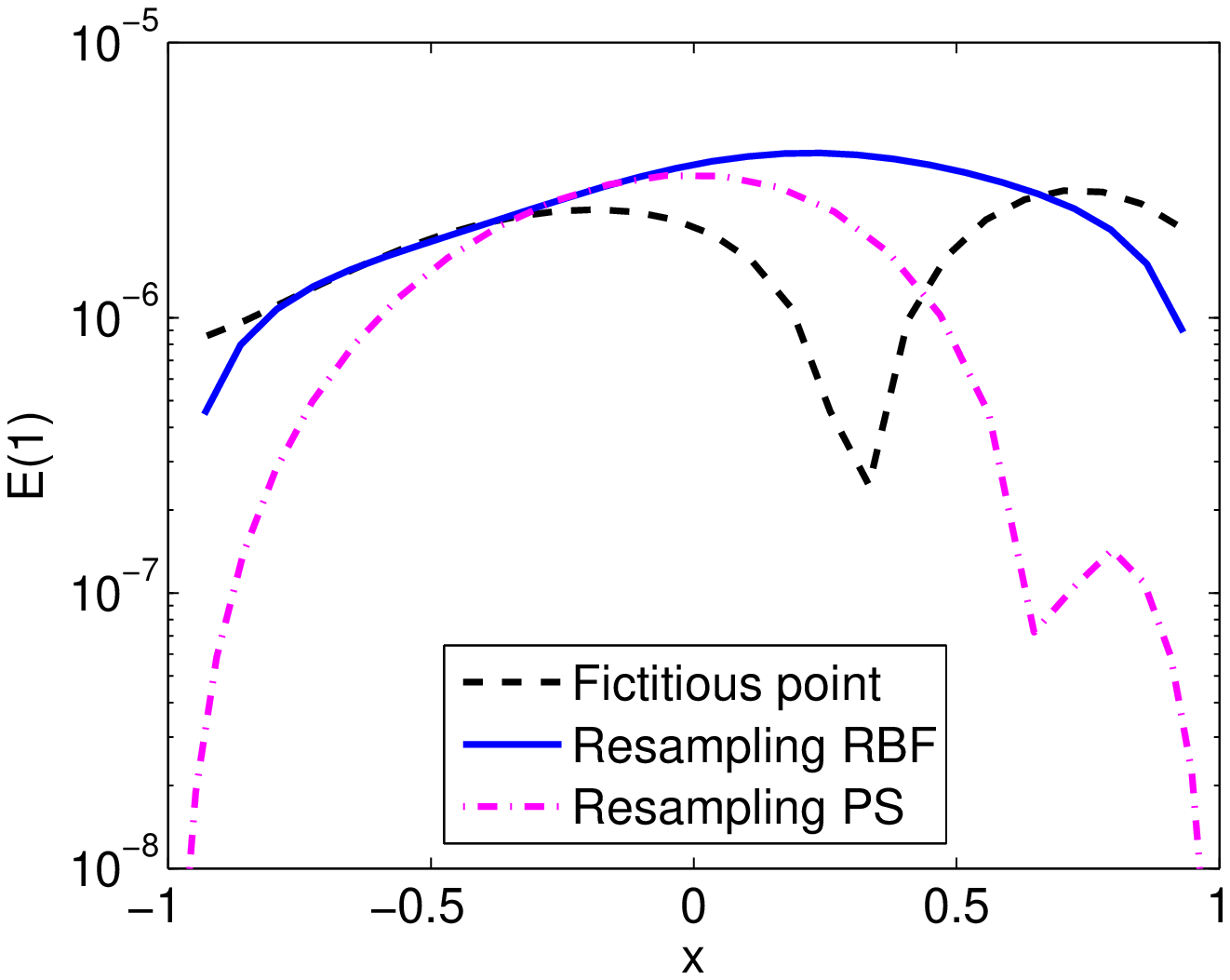}
\includegraphics[width=0.49\textwidth]{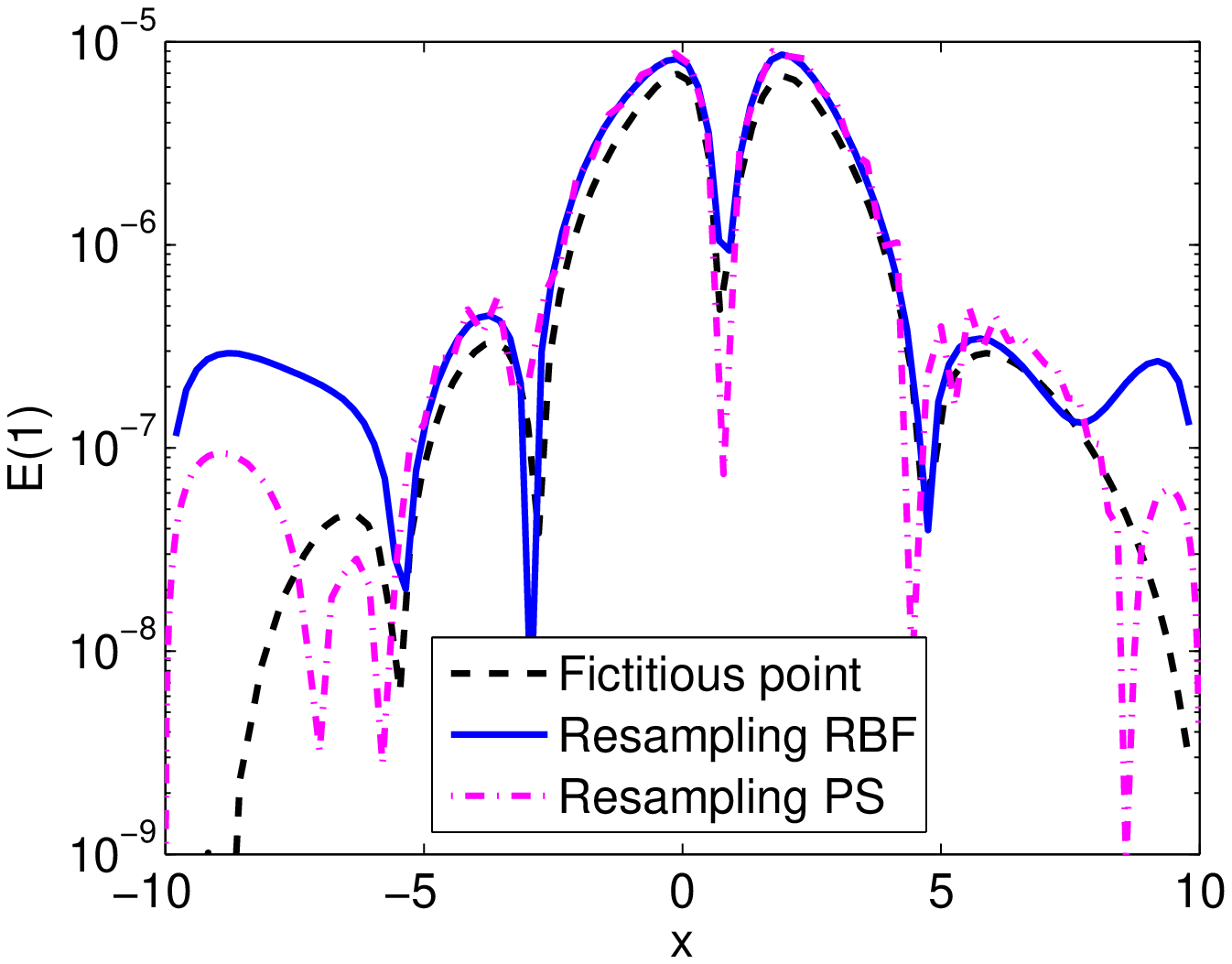}
\caption{Absolute error of the fictitious point RBF method, the resampling RBF method, and the resampling pseudo-spectral method at time $t=1$ for $L=1$ and $N=30$ (left) and for $L=10$, $N=100$ (right). For the RBF methods $\ep{}=\frac{0.08}{h}$ was used.}
\label{fig:absolute}
\end{figure}

To illustrate the capability of the proposed methods, we start with comparing the approximation accuracy with that of a pseudo-spectral resampling method. The pseudo-spectral method employs the same number of Chebyshev nodes as the number of uniform nodes used by the RBF methods. For a description of its implementation, see~\cite{Drishal}.
The absolute errors for two different values of $L$ are plotted in Figure~(\ref{fig:absolute}). As shown in the figure, the errors of both the RBF based methods and the pseudospectral resampling method are similar in magnitude. For the shorter interval, the pseudo-spectral method has smaller errors near the boundaries, which is consistent with the clustering of the Chebyshev nodes. However, for the larger interval, where the solution is small at the boundary, this effect is not visible.

\begin{figure}[h!]
\centering
\includegraphics[width=0.49\textwidth]{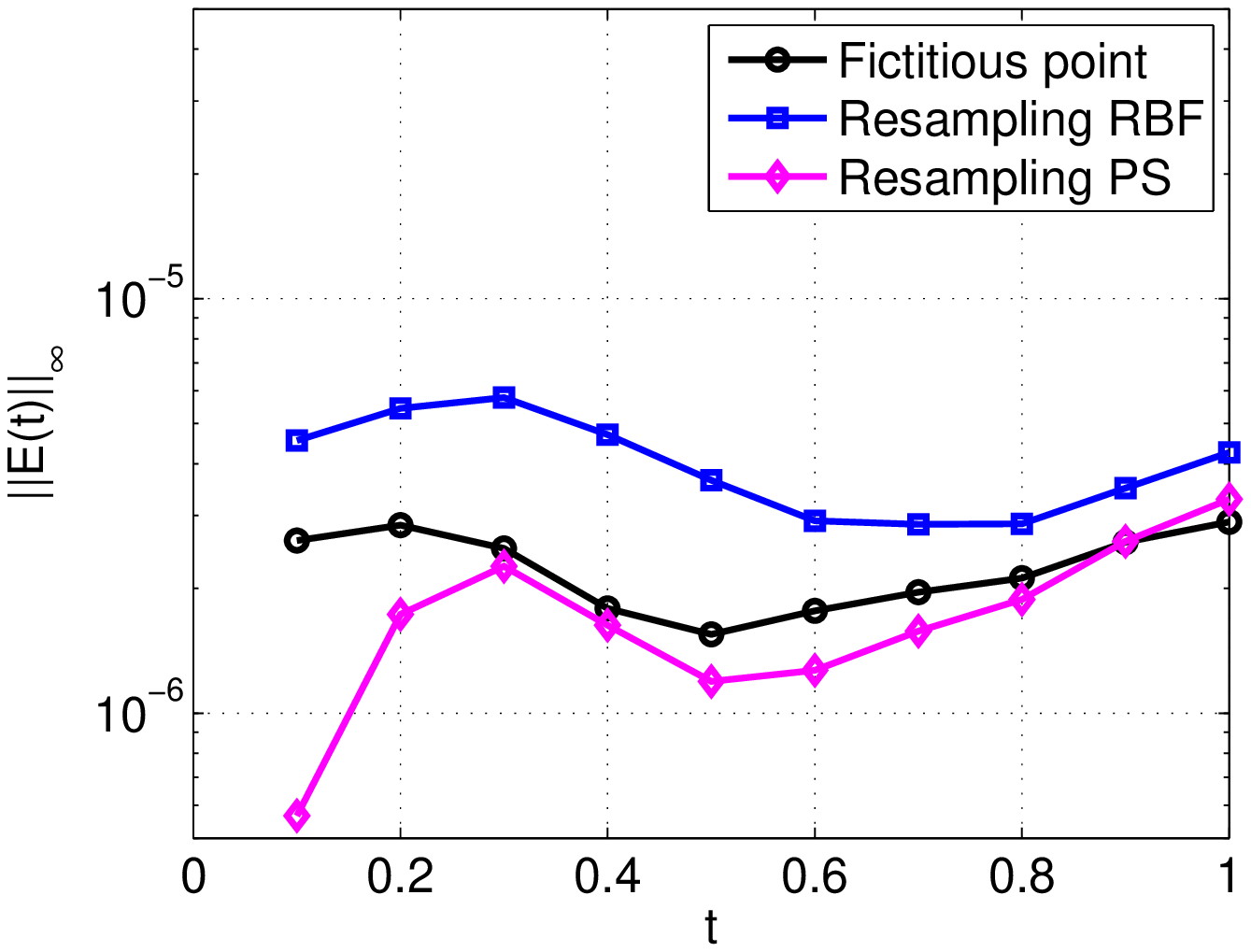}
\includegraphics[width=0.49\textwidth]{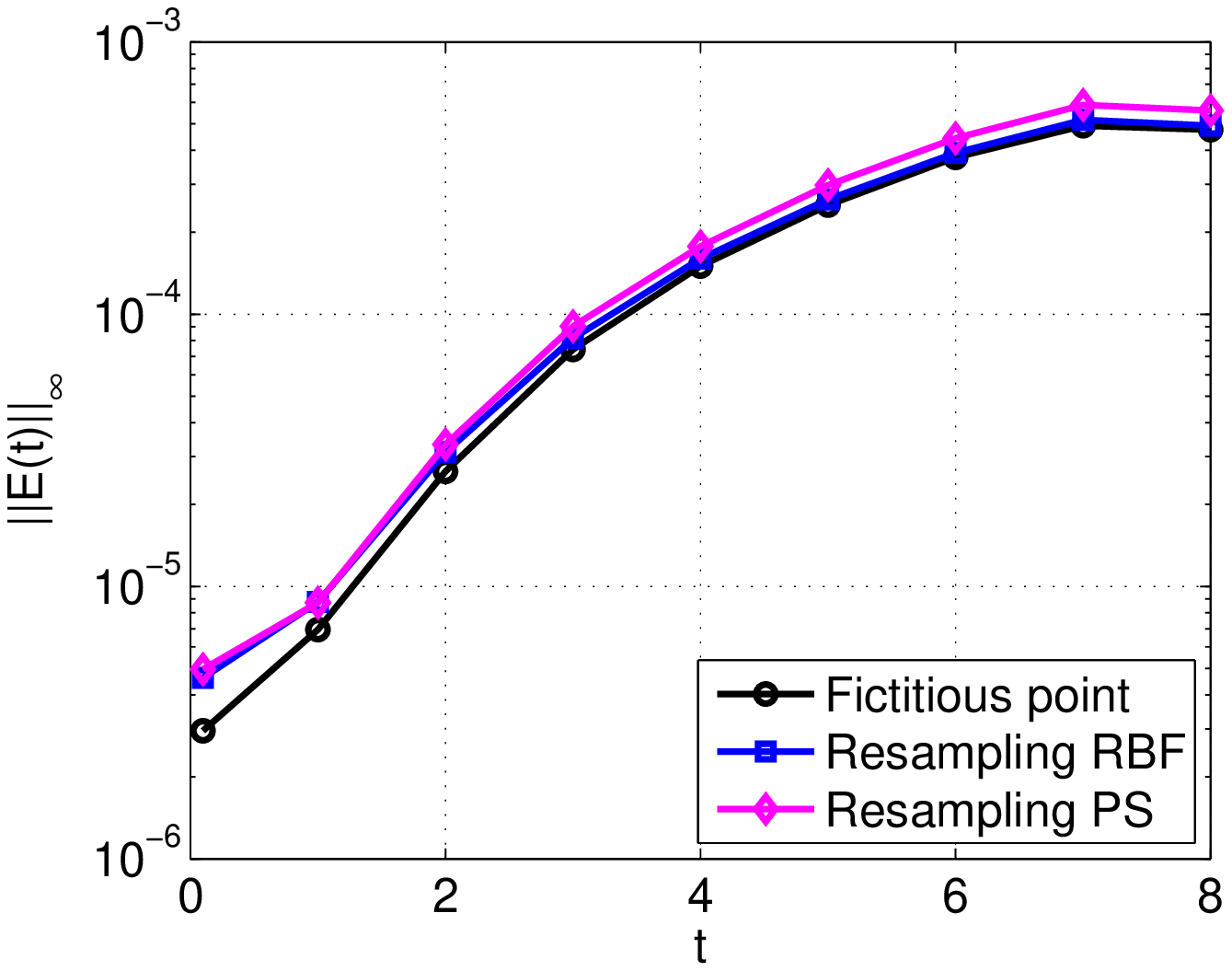}
\caption{The $L_{\infty}$ error as a function of time $t$ for the fictitious point RBF method, the resampling RBF method, and a resampling PS method. Results are shown for $L=1$ and $N=30$ (left) and $L=10$ and $N=100$ (right). Both RBF methods use $\ep{}=\frac{0.08}{h}$ and a uniform node distribution, while the PS method employs a Chebyshev node distribution.}
\label{fig:errort}
\end{figure}

The $L_{\infty}$ errors over time for the approximated solutions are illustrated in Figure~\ref{fig:errort}.
%Uniform nodes, $N=100$, are used to discretize the domain with $\ep{}=0.5$ for fictitious
%point method and resampling RBF method respectively. In the resampling PS method, $N=100$ Chebyshev nodes is used to produce the approximation solution. In this figure, the domain of computation is $\Omega=[-L,L]$.
% g(u) contains u^5 which gives q=4, q(1-q/5)=0.8
% exp(c_3 (1+t).../...)=exp(c_3 (1+t)^1.8/1.8) c3=tilde(cq)*Qinv*max(Psi,B) Cq=24*C Ct=(q+1) 120*2+max(F) Qinv*Psi=100
If we go back to the global error estimate~\eqref{eq:gdiserror}, and insert $q=4$ (for this test case), the exponential time-dependent growth factor becomes $\exp(C_3((1+t)^{1.12}-1)/1.12)$. We do not know the precise value of $C_3$, but based on our experiments a value larger than one should be expected, in which case the predicted growth would be at least two orders of magnitudes larger than what is observed. However, the growth factor in the estimate arises from the growth estimate~\eqref{eq:uxt} for the solution $u(x,t)$, which in our particular case does not grow at all. Hence, it would be possible to make a more  favorable estimate for our special case. Both the accuracy and the growth rate of the errors of the three different solutions are similar. For the shorter interval $L=1$, the resampling RBF method is slightly worse than the other two methods, while for the longer interval $L=10$, both RBF methods are slightly more accurate than the pseudo-spectral method.
% -1/h goes crazy because h is defined differently

\begin{figure}[h!]
\centering
\includegraphics[width=0.49\textwidth]{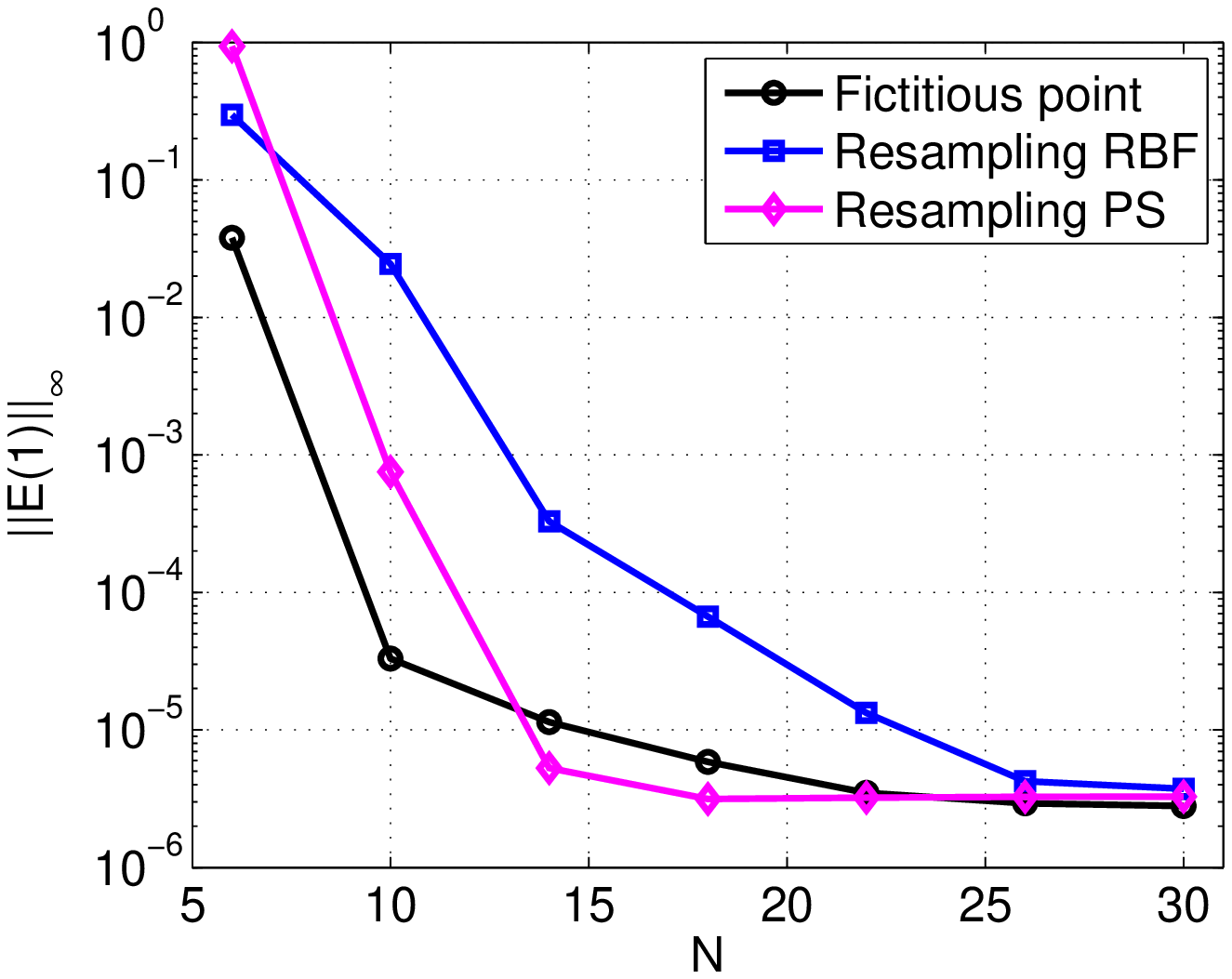}
\includegraphics[width=0.49\textwidth]{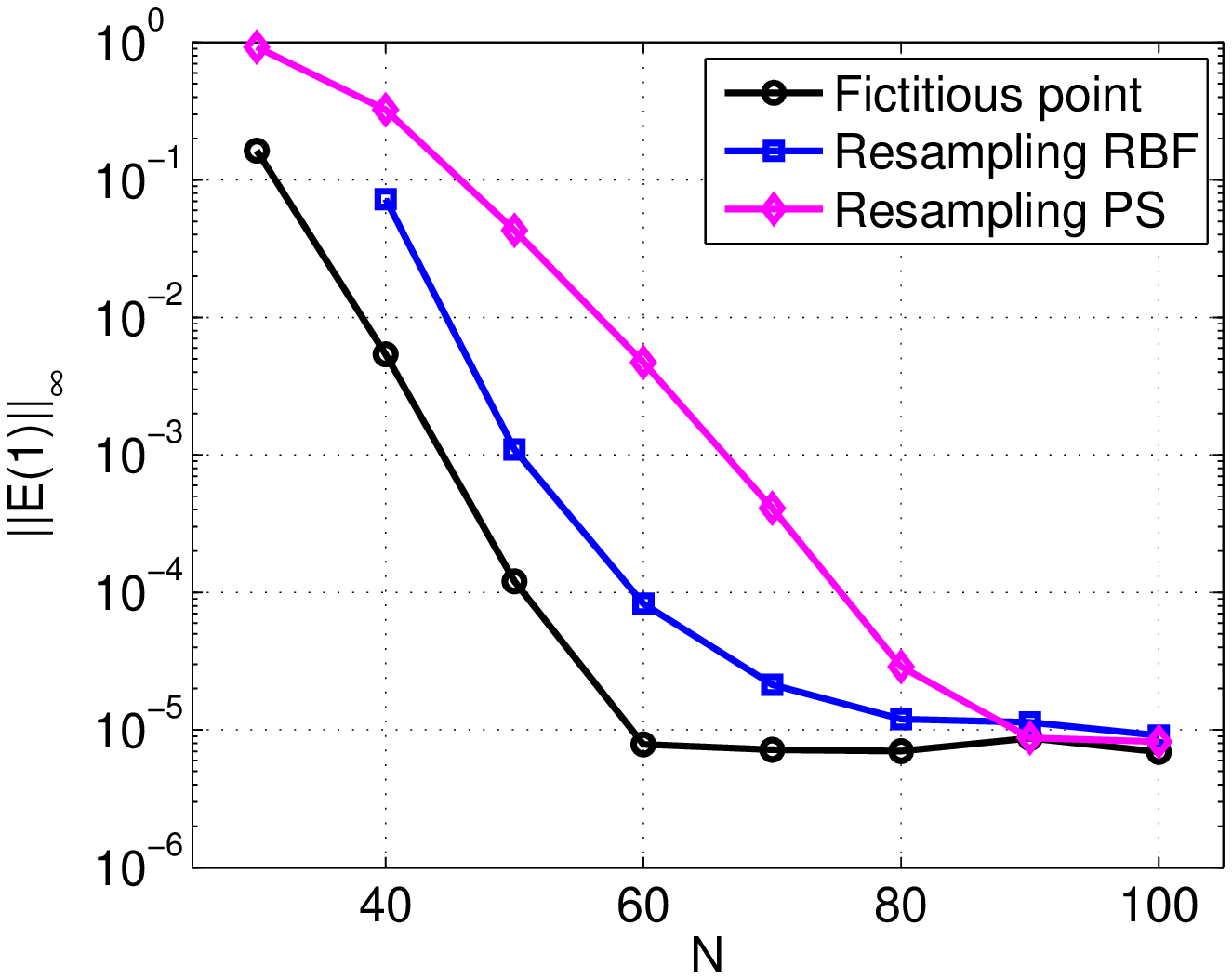}
\caption{The $L_{\infty}$ error at time $t=1$ as a function of the number of node points $N$ for the fictitious point RBF method, the resampling RBF method, and a resampling PS method. Results are shown for $L=1$ (left) and $L=10$ (right). Both RBF methods use $\ep{}=\frac{0.08}{h}$ and uniform node distributions, while the PS method employs Chebyshev node distributions.
}
\label{fig:conver}
\end{figure}

Figure~\ref{fig:conver} displays the convergence behavior as a function of $N$ for the two RBF methods compared with the PS method. For all three methods, the highest attainable accuracy is the same. When $\ep{}h$ is constant, as in this experiment, we would expect the error to reach a saturation level, but accuracy is also limited by conditioning, and this may be the effect that we see here. In both cases, the fictitious point RBF method reaches the highest accuracy faster than the  resampling RBF method. The PS method performs best for the short interval, and performs worse than the RBF methods for the longer interval. One explanation for this can be that the node density for the Chebyshev nodes compared with the uniform nodes is lower in the interesting region (middle of the domain) in this case.
%
%In Resampling PS method for conditioning of mass matrix in ode solver, we can not show the convergence rate for a big number of Chebyshev nodes.
%Scaling the shape parameter in Resampling RBF method and
%Fictitious point method allow us to have a convergence rate for a
%different data sets. In Figure~(\ref{fig:conver}), from left to
%right, we use the uniform distribution of nodes for Resampling RBF
%and Fictitious point methods with $\varepsilon=1$ and
%$\varepsilon=1.5$. From Figure~(\ref{fig:conver}), it can be seen
%that the $L_{\infty}$ error decrease as the number of collocation
%points increases. Figure~(\ref{fig:epsi}) shows the $L_{\infty}$
%error as a function of $\varepsilon$ for different uniform
%distribution of data sets.
%

\subsection{Two-dimensional square domain}
\begin{figure}[h!]
%\centerline{
\centering
\includegraphics[width=0.49\textwidth]{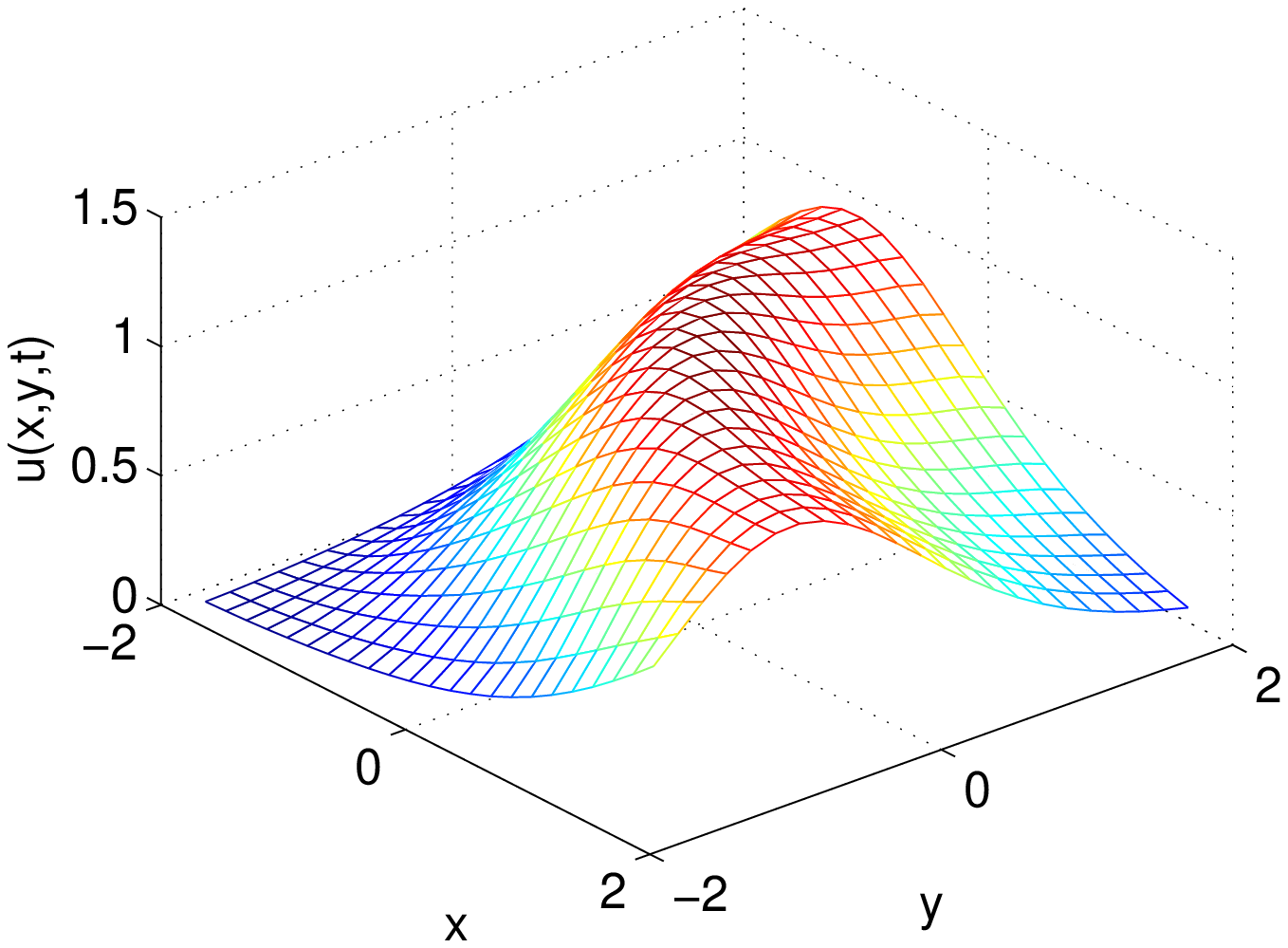}
\includegraphics[width=0.49\textwidth]{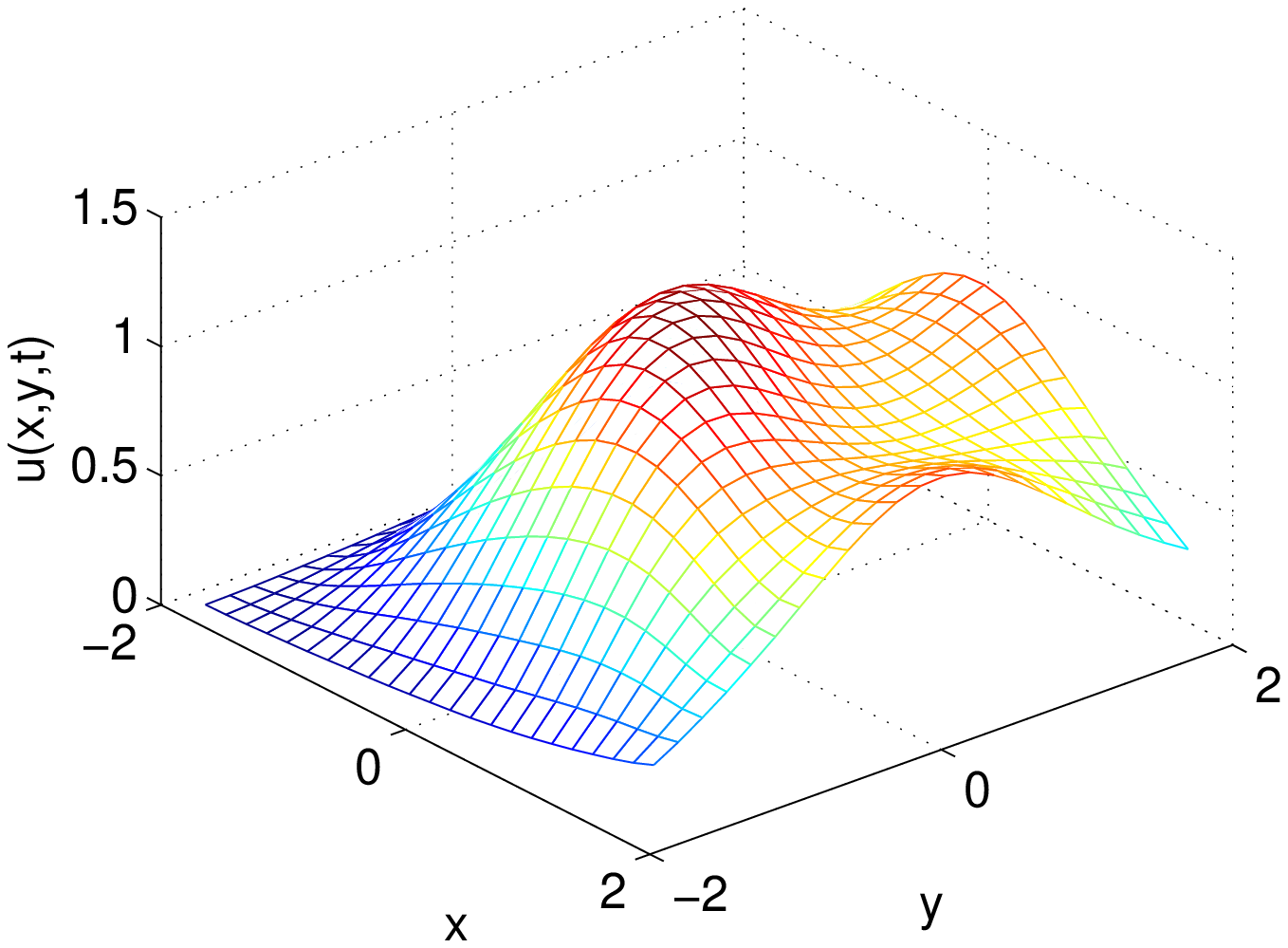}
\caption{Resampling RBF approximate solution in the square domain $\Omega=[-2,\,2]^2$ at time $t=1$ (left) and $t=2$ (right) with $N=25\times 25$ points and shape parameter $\ep{}=0.8$.}
\label{fig:2dapp}
\end{figure}

% Talk about the flexibility of the RBF
In this section, we demonstrate how the flexibility of the RBF approximations allows us to implement the fictitious point method and resampling RBF method in a two-dimensional domain with a similar effort as for the one-dimensional problem. Here, we do not compare with the resampling PS method, which is less straightforward to implement. We consider the square domain $\Omega=[-L,\, L]\times [-L,\, L]$ and the Rosenau equation \eqref{eq:Rosenau} with $\alpha=1$, $g(u)=u^3+u^2$ and initial condition $f_0(x,y)=\textrm{sech} (x+y)$ and boundary conditions
\begin{align}
f_1(x,y,t)&=\textrm{sech}(x+y-t),\\
f_2(x,y,t)&=-\textrm{sech}(x+y-t)\textrm{tanh}(x+y-t),
\end{align}
where the derivative in the second condition is taken as either $u_x$ or $u_y$ depending on which part of the boundary is involved.
%In this case, an analytical solution is not available.
For the two-dimensional test cases, we do not have any analytical solutions. The approximate solution at two different times is displayed in Figure~\ref{fig:2dapp}.

\begin{figure}[h!]
%\centerline{
%\includegraphics[width=0.6\textwidth]{squarfic2.eps}
%\includegraphics[width=0.6\textwidth]{squarres2.eps}}
\centering
\includegraphics[width=0.4\textwidth]{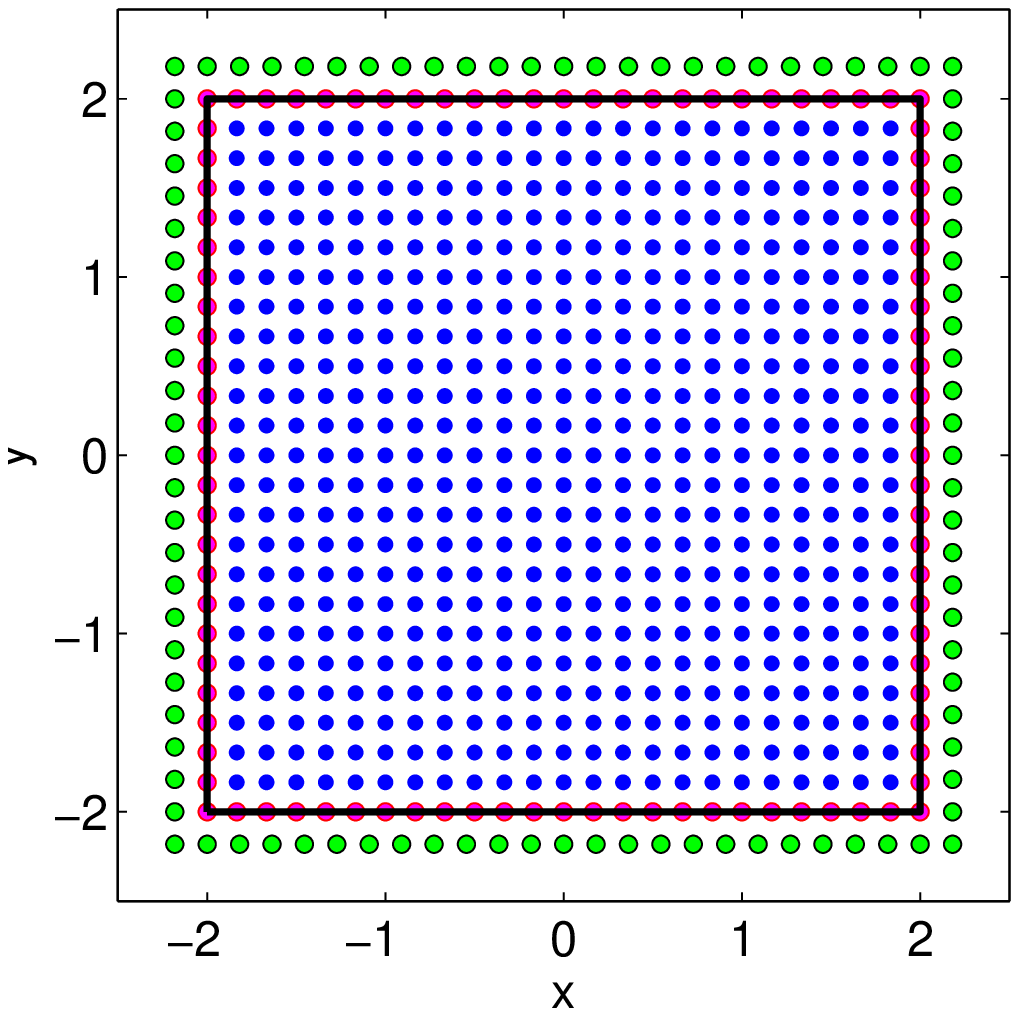}
\includegraphics[width=0.4\textwidth]{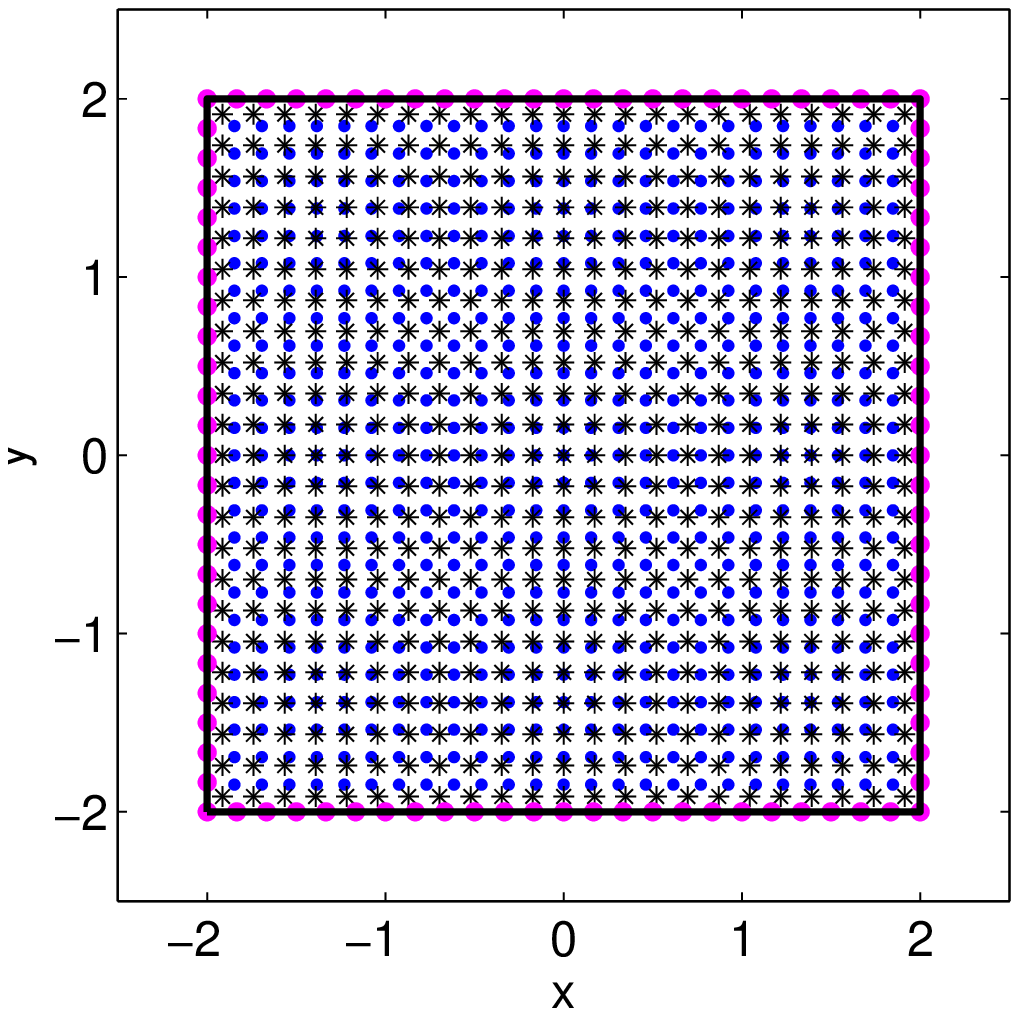}
\caption{A node distribution for the fictitious point method, where the fictitious points are uniformly distributed outside the domain (left) and a node distribution for the resampling RBF method with auxiliary points displayed as $\ast$ (right) for a square computational domain. }
\label{fig:squardom}
\end{figure}

We start from a uniform discretization of the domain $\Omega$ with $N=n^2$ points. We denote the number of interior points by $N_d$ and the number of boundary points by $N_b$. 
For the fictitious point method we need to add $N_b$ fictitious points outside the domain to enforce the Neumann boundary conditions. Note that if we simply choose an extension of the uniform grid, the resulting number of fictitious points is too large. For the resampling RBF method we generate $N_d-2N_b$ auxiliary points inside of the domain. Note that here the number of boundary points is modified (and they are therefore not on the uniform grid) to make the interior and auxiliary node numbers compatible. Sample node distributions for the two methods are displayed in Figure~\ref{fig:squardom}.

\begin{figure}[h!]
\centering
\includegraphics[width=0.49\textwidth]{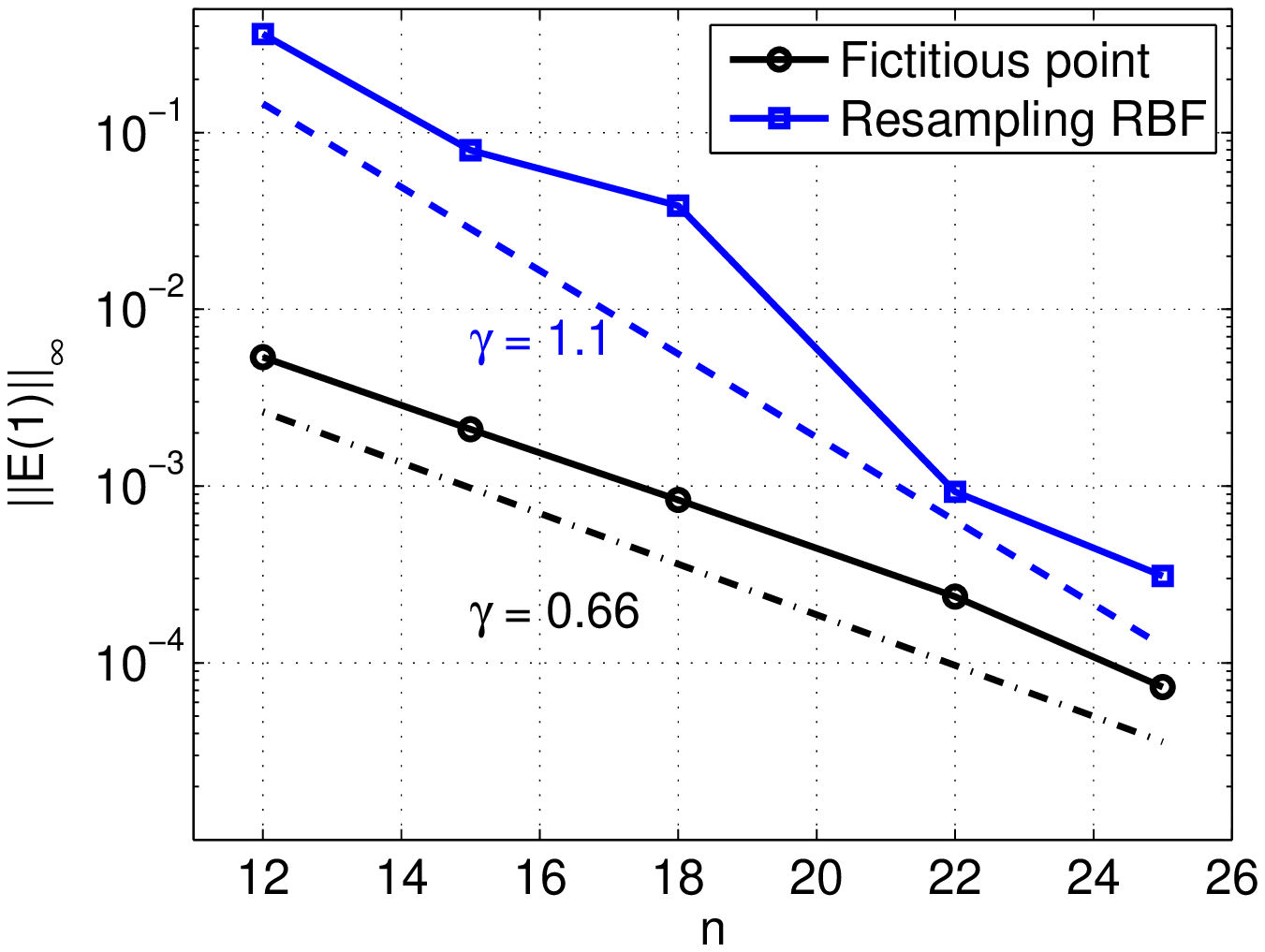}
\includegraphics[width=0.49\textwidth]{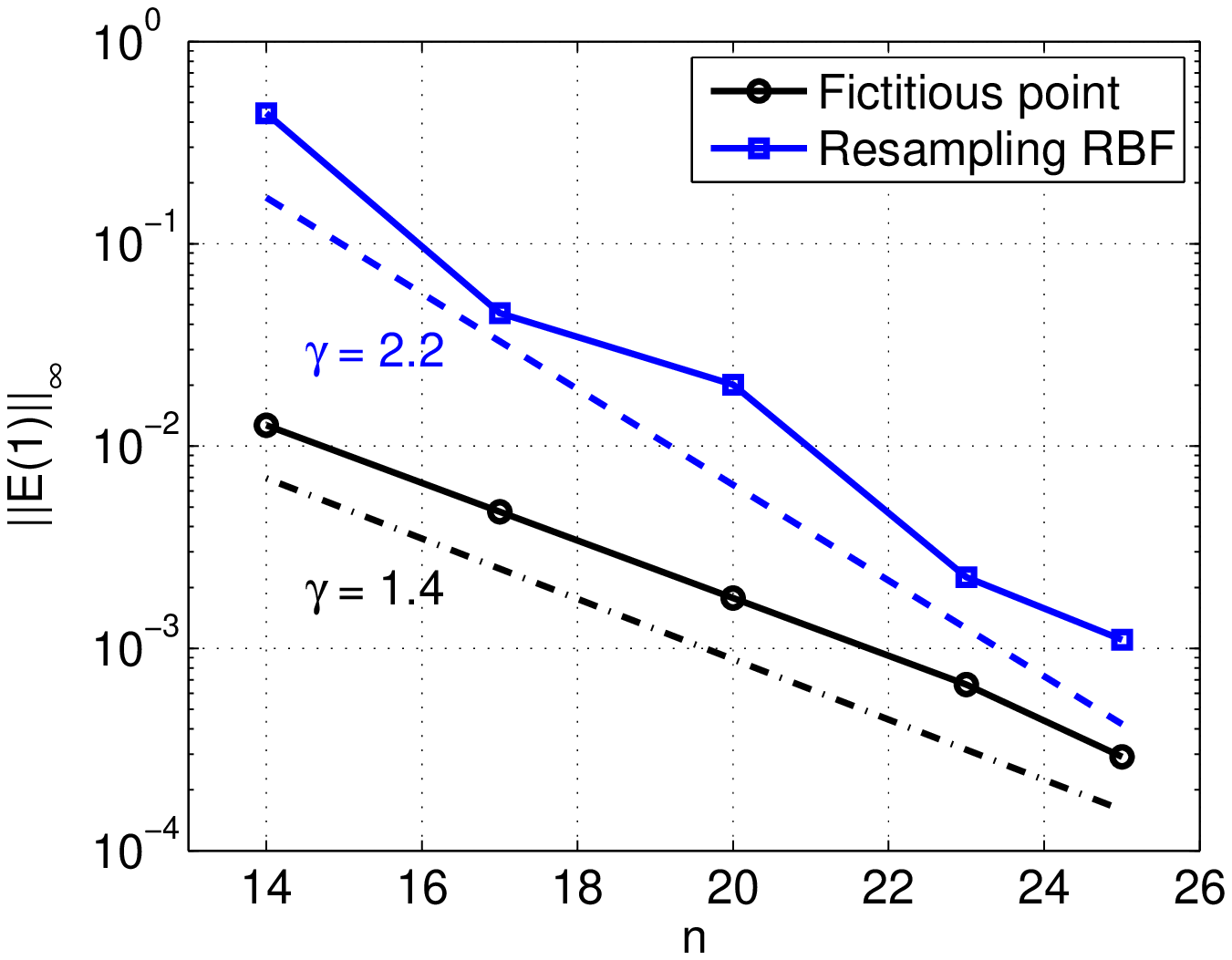}
\caption{The error at time $t=1$ against the number of points per dimension $n$ for $L=1$ and $\ep{}=1.6$ (right) and $L=2$ and $\ep{}=0.8$ (left). In both cases, errors are computed against a reference solution computed with the fictitious point RBF method for $n=28$, and the error is evaluated at $25\times25$ interior points.}
\label{fig:conver2D}
\end{figure}

According to the error estimate~\eqref{eq:gdiserror} for the fictitious point method, we expect exponential convergence in $1/\sqrt{h}$ for fixed $\ep{}$. In practice, we often observe exponential convergence in $1/h$. In  Figure~\ref{fig:conver2D}, we plot the error as a function of $n\propto 1/h$. For this range of $n$-values, the conditioning is low enough to not influence the error, and exponential convergence can be observed for both RBF methods. We see that the estimated slopes in the right subfigure are precisely double those in the left subfigures. If we take into account that $h$ is also twice as large for the case $L=2$, we can conclude that the rate of convergence in terms of $1/h$ is the same in both cases. Again, the fictitious point is more accurate, while the resampling method converges faster. The fictitious point uses a larger number of basis functions, which could account for the smaller error, but it is not clear why the rates differ in the way they do.
%The error is significantly smaller for the fictitious point method for lower $n$, while the convergence rate is higher for the resampling method such that both are equally accurate for $n=25$. Hence, which method to chose depends on the resolution of the problem.
%shows the convergence trends versos number of points for the fixed value of shape parameter. The max errors are computed in both cases based on equally spaced points. The resampling RBF with $n=35$  is used to produce the reference solution.

%From the figure, it can be seen that in both methods we have almost the same accuracy for large number of the points but we can expect the good accuracy from the fictitious point method for the lower number of the points in compare the resampling RBF method.
%Note: As we see from the figure like the one dimensional case the fictitious point method is stable when we keep the $\ep{}$ fixed and change the number of the points in large range (in other word, the good accuracy can be achieved for large range of $\ep{}$), but for the resampling method we can not keep $\ep{}$ fixed for the large range of number of points, and for having the good accuracy need to use proper $\ep{}$ for related the $N$.

%( We can have the figure od convergence against h or N. Please choose every one that is better.)

%\begin{figure}[h!]
%\centerline{
%\includegraphics[width=0.5\textwidth]{conver2dh1.eps}
%\includegraphics[width=0.5\textwidth]{conver2dh2.eps}}
%\caption{Error as function of $h$ where $L=1$, $t=1$ and  from left to right $\ep{}=1.5$, $\ep{}=2$. In both cases, errors are computed by resampling RBF reference solution of the corresponding type with $n=35$.}
%\label{fig:conver2D}
%\end{figure}

\subsection{Two-dimensional starfish domain}
\begin{figure}[h!]
%\centerline{
%\includegraphics[width=0.35\textwidth]{starfishfict05.eps}
%\includegraphics[width=0.35\textwidth]{starfishfict1.eps}
%\includegraphics[width=0.35\textwidth]{starfishfict2.eps}}
\includegraphics[width=0.33\textwidth]{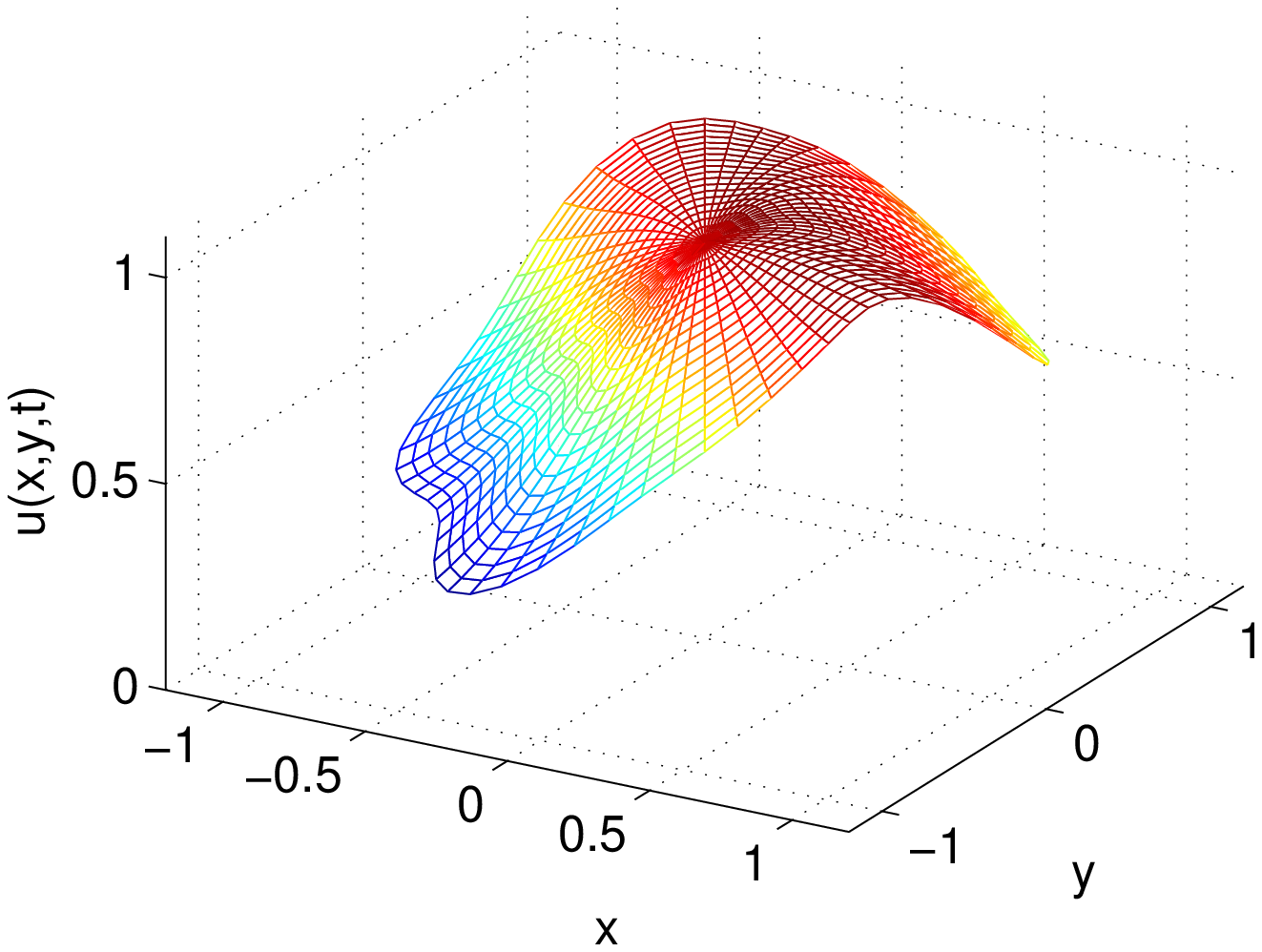}
\includegraphics[width=0.33\textwidth]{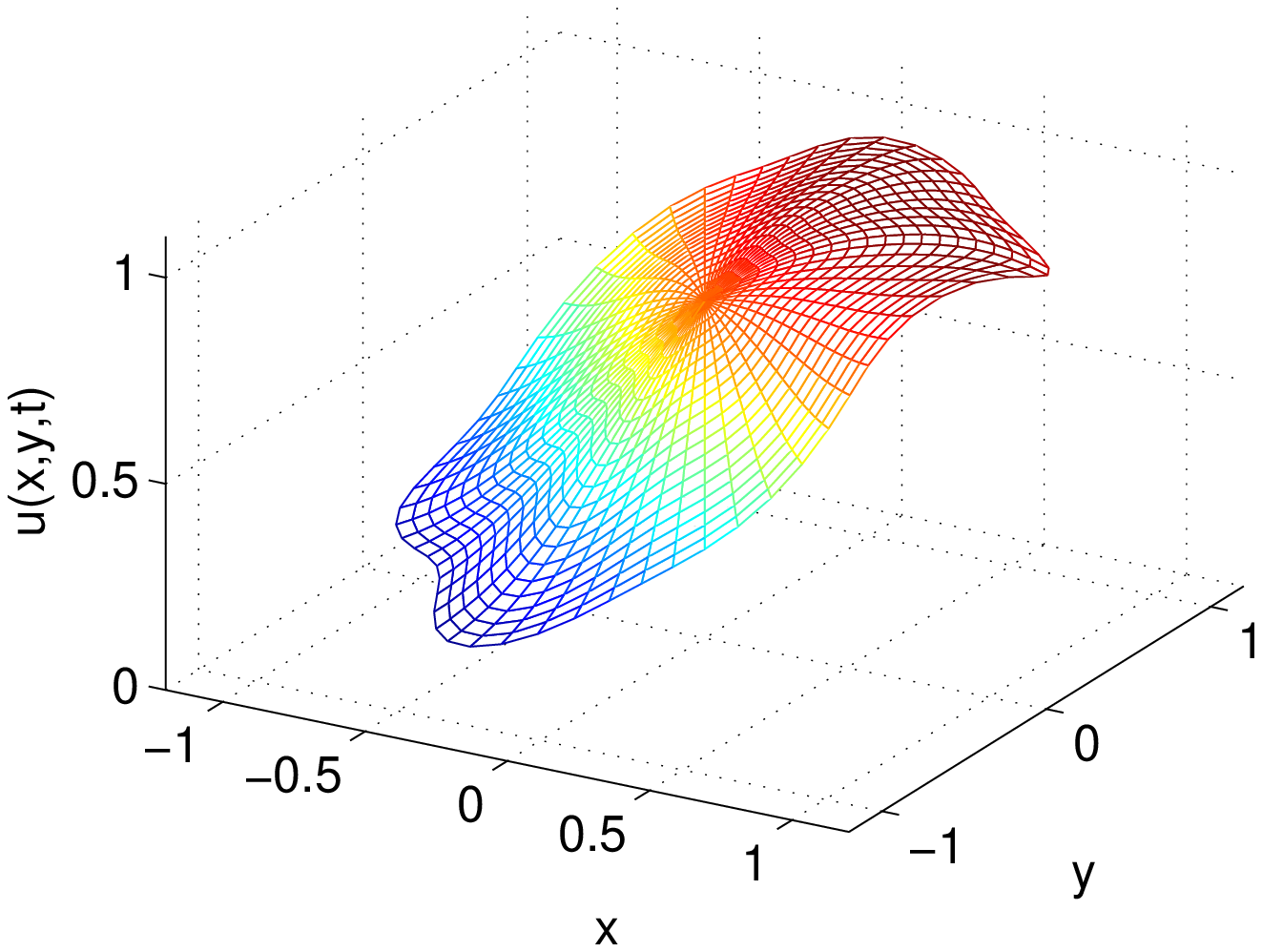}
\includegraphics[width=0.33\textwidth]{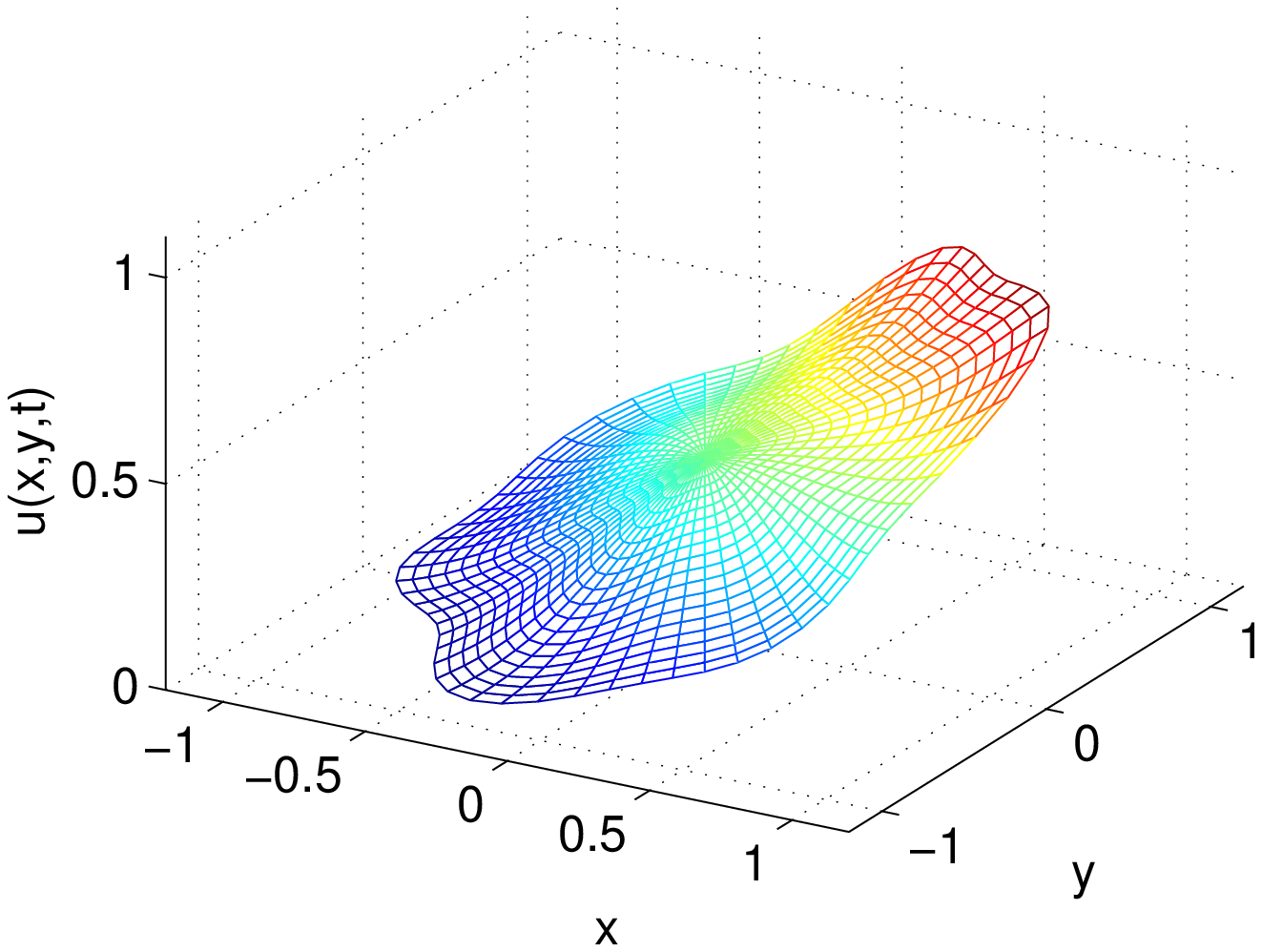}
\caption{The approximate solution in the starfish domain with $N_d=253$ interior points, $N_b=56$ boundary  points using $\ep{}=1.5$ for
%that is shown in left  figure~\ref{fig:starfishdom} where $\ep{}=1.5$. 
$t=0.5$ (left), $t=1$ (middle), and $t=2$ (right).}
\label{fig:starfishsol}
\end{figure}

% grad u = (1,1)/sqrt(2)*-sech(x+y-t)*tanh(x+y-t)
We now take a step further in demonstrating the flexibility of the RBF based methods by considering an irregular two-dimensional domain.
As a test problem, we consider the starfish domain with boundary defined by the parametric equation
\begin{equation}
%r(\theta)=1+0.1(\sin(6\theta)+\sin(3\theta)),\quad \theta\in[0,2\pi),
r(\theta)=1+0.07(\sin(6\theta)+\sin(3\theta)),\quad \theta\in[0,2\pi).
\end{equation}
We also need the derivative of the boundary equation in order to compute the outward normal direction $n=(n_x,n_y)$, which is needed for the boundary conditions. We have
\begin{equation}
(n_x,n_y) = \frac{r^\prime(\theta)(\sin(\theta),-\cos(\theta)) + r(\theta)(\cos(\theta),\sin(\theta))}{\sqrt{r^\prime(\theta)^2+r(\theta)^2}}.
%(r'sin t + rcos t)^2+(rsint-r'cost)^2=r'^2+r^2 2r'rsintcost-2r'rsintcost 
%r(\theta)=1+0.1(\sin(6\theta)+\sin(3\theta)),\quad \theta\in[0,2\pi),
\end{equation}
We use a similar test problem as for the square domain with boundary Dirichlet data 
\begin{equation}
f_1(x,y,t)=\mathrm{sech}(x+y-t). 
\end{equation}
For the normal derivative condition, we impose
%for the Rosenau equation~\eqref{eq:Rosenau} with $\alpha=1$, $g(u)=u^3+u^2$ accompanied by initial condition $f_0(x,y)=\textrm{sech} (x+y)$ and boundary conditions
\begin{equation}
%f_1(x,y,t)&=\textrm{sech}(x+y-t), & (x,y)&\in\partial\Omega,\\
f_2(x,y,t)=\nabla u\cdot n=-\textrm{sech}(x+y-t)\textrm{tanh}(x+y-t)(n_x+n_y).
\end{equation}
The approximate solution at three different times is shown in Figure~\ref{fig:starfishsol}.
%%%%%%%%%%%%%%%%%%%%%%%%% HERE

\begin{figure}[h!]
\centering
\includegraphics[width=0.49\textwidth]{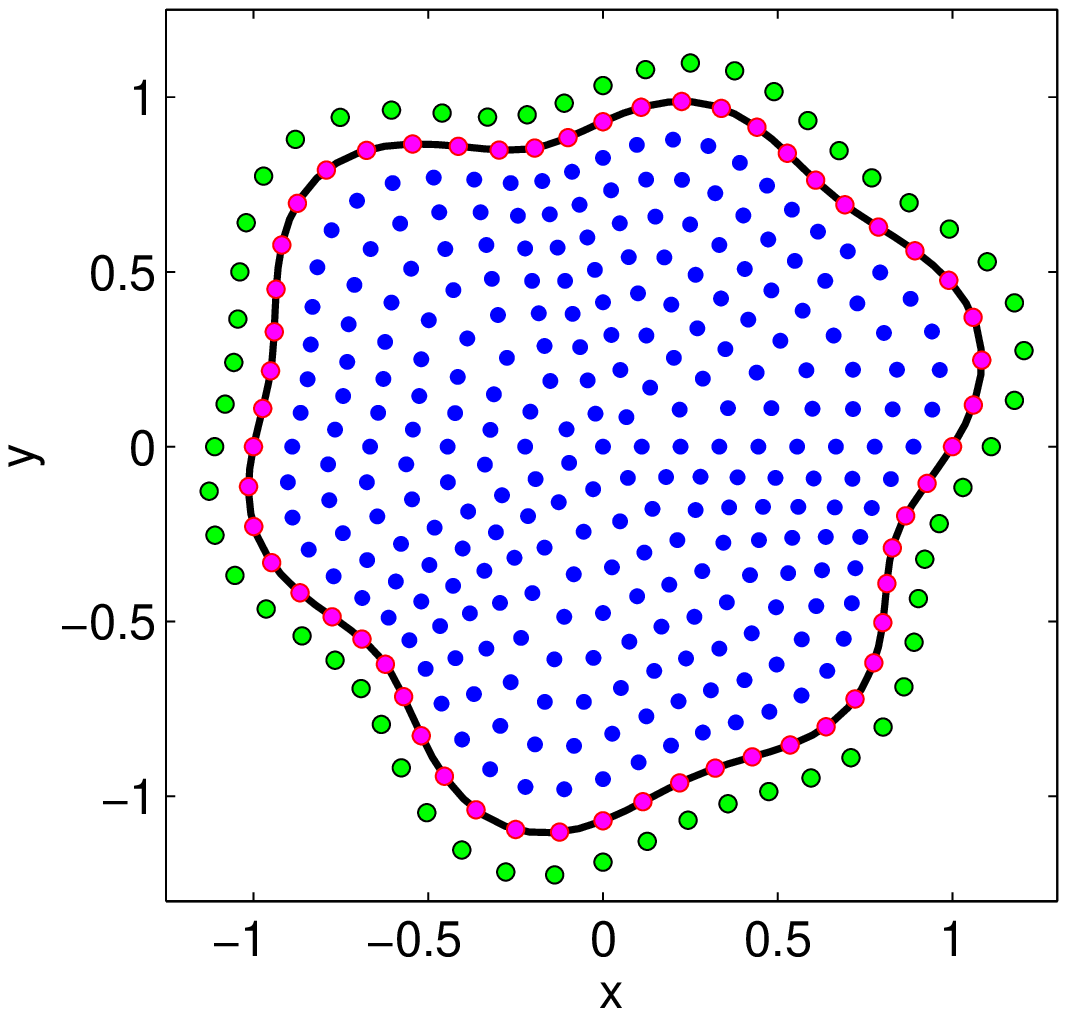}
\includegraphics[width=0.49\textwidth]{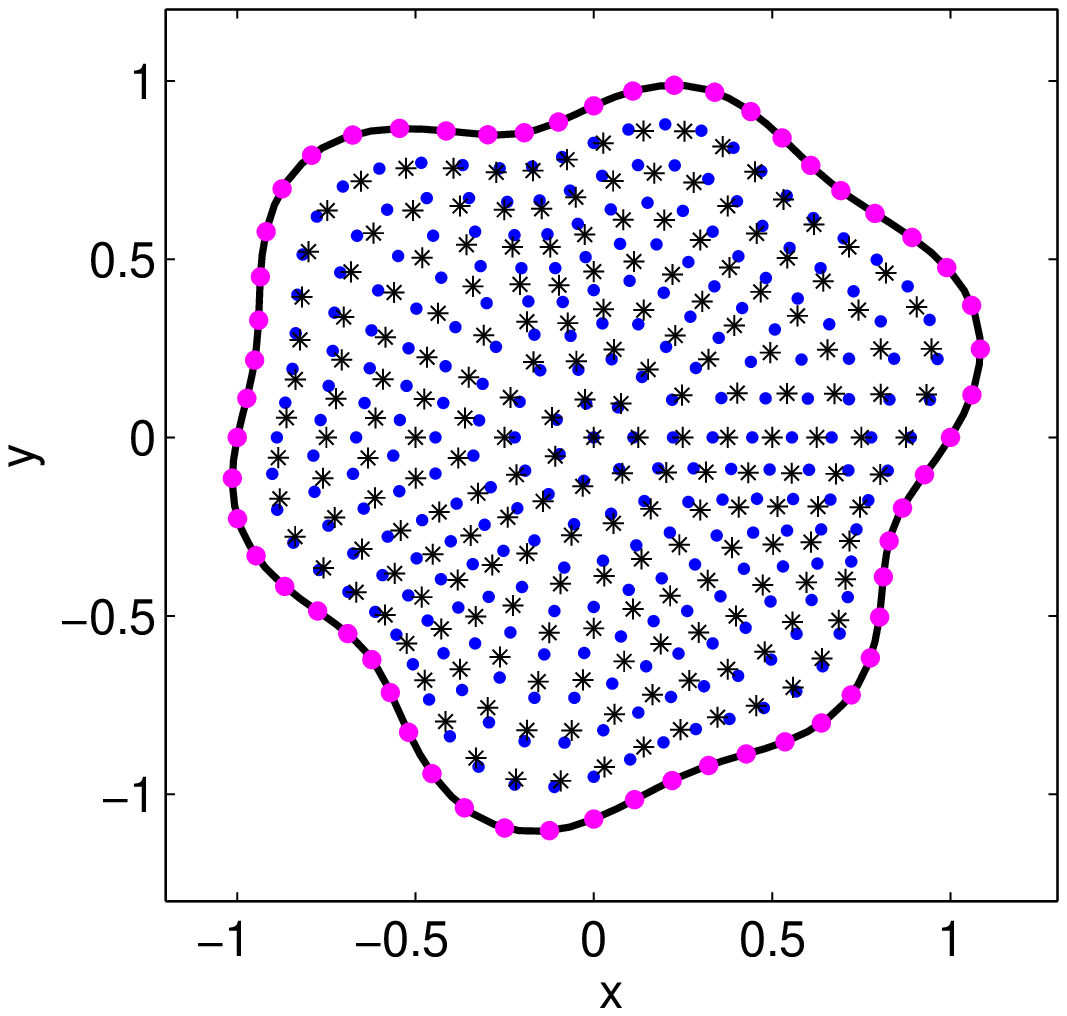}
\caption{Sample node distributions with $N_d=253$ and $N_b=56$ for the starfish domain for the fictitious point method (left) and the resampling RBF method (right) with auxiliary points marked with $\ast$.}
\label{fig:starfishdom}
\end{figure}

Sample node distributions for the fictitious point and resampling RBF methods are illustrated in Figure~\ref{fig:starfishdom}.
Just as for the square domain, the number of fictitious points outside the domain is $N_b$ and the number of auxiliary points inside the domain in the resampling method is $N_d-2N_b$.

\begin{figure}[h!]
\centering
\includegraphics[width=0.49\textwidth]{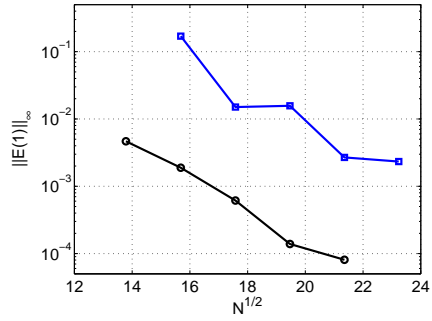}
\caption{Error as a function of the square root of the number of points $\sqrt{N}=\sqrt{N_d+N_b}$ where $t=1$ and $\ep{}=2$. The reference solution is produced using the fictitious point method with $N=540\approx 23.24^2$.}
\label{fig:converstar}
\end{figure}

The max error as function of $N$ for a fixed shape parameter value, $\ep{}=2$, is illustrated in Figure~\ref{fig:converstar}. The reference solution is computed using the fictitious point method with $N=540$ nodes. The max error is estimated from evaluation at $600$ radially uniformly distributed points in the domain.
The error trends are similar to those for the square domain, showing that the RBF methods provide a well functioning generalization of both the fictitious point method and the resampling method to general domains.

\section{Conclusion}
\label{sec:conc}
% Multiple boundary conditions provide an extra difficulty, several conditions
The Rosenau equation, which is used as an application throughout this paper, is an example of a non-linear PDE with multiple boundary conditions as well as mixed space-time derivatives. Multiple boundary conditions provides an extra challenge when solving PDE problems. The standard form of a typical collocation method assumes that one condition is imposed at each node point/grid point. Hence, the additional conditions at the boundary nodes lead to a mismatch between the number of conditions and the number of unknowns. 
% Where does the rosenau come in? Example used throughout the paper

% Introduced for PS-methods.
Two approaches to manage multiple boundary conditions that have introduced for spectral methods are fictitious point methods and resampling methods. In this paper we have shown how to implement these approaches in the context of RBF collocation methods.
% Here we have shown how to implement these with an RBF collocation method.
% In 1-d accuracy similar
From numerical experiments for a one-dimensional test problem, we could see that the behavior of the method with respect to accuracy in space and time is very similar to that of the corresponding pseudo-spectral method.
% In 2-D the matching of conditions becomes difficult. We can do it.

For two-dimensional problems, already in a regular geometry such as the square, the application of spectral methods becomes more complicated. Approximations are typically based on tensor product grids, but if we use the one-dimensional extension techniques for each grid line, again the number of extra conditions and extra points do not naturally match. The problem can for example be resolved by choosing one of the directions for the corner points, but then the approximations in the other direction needs to be of lower order. 

% Require regular geometry. We can do irregular.
We show that with the two RBF methods, due to the freedom of node placement, we can distribute the fictitious points or the resampled nodes uniformly and symmetrically with respect to the domain. Furthermore, we show that the concept can be transferred also to irregularly shaped domains.

% Theoretical results given estimates for u, we can also estimate the error.
We have also analyzed the theoretical properties of the fictitious point RBF approximation for the one-dimensional Rosenau equation. We could show that the spectral convergence of the spatial approximation carries over to the PDE solution, while the growth of the error in time in our estimate strongly depends on the growth estimates for the exact solution of the problem. Using a generally valid growth estimate, as we did, led to an overestimate of the error growth in time for the test problem we used, where the norm of the exact solution does not grow. 

To conclude, both the implemented approaches are promising for problems with multiple boundary conditions, especially for geometries where spectral methods cannot easily be applied. Global RBF approximations as the ones used here are competitive for problems in one or two space dimensions, but the computational cost can become prohibitive for higher-dimensional problems due to the need to solve dense linear systems. Therefore, an interesting future direction is to see how resampling and fictitious point methods can be combined with localized (stencil or partition based) RBF methods.

%The fictitious point and resampling RBF is introduced as efficient methods for implementation of multiple boundary conditions and applied to the Rosenau equation for a test problem. A combination of theoretical and experimental analysis is studied for the fictitious point method. The estimated error indicates that we can achieve high order of convergence rate which can be slower by time.

%For one-dimensional case,  the error comparison with resampling PS method verifies the accuracy of the proposed RBF based method. The numerical evidence ensure us the acceptable accuracy with a few number of points can be achieved by fictitious point and resampling RBF method.

%Easy implementation for one dimensional case  as well as two dimensional case and flexibility of implementation in the irregular domains such as starfish like domain with uniform and non-uniform disretization can be counted as advantage of the fictitious point and resampling RBF approximation method for numerical investigation of the Rosenau equation.
%
%%%%%%%%%%%%%%%%%%%%%%%%%%%%%%%%%%%%%%%%%%%%%%%%%%%%%%%%%%%%%%%%%%%%%%%%%%%%%%%%%%%%%%%%%%%%%%%%%%%%%%%%
%\bibliographystyle{spmpsci}
%\bibliography{refs}

\end{document}